\input ./style/arxiv-general.cfg
\documentclass[aap,MSNbibl,seceqn,dvips]{arximspdf}
\makeatletter
   \@ifpackageloaded{graphicx}{}{\usepackage{graphicx}}
\makeatother
\usepackage{mathbh}

\doi{10.1214/14-AAP1036} 
\volume{25}
\issue{4}
\pubyear{2015}
\firstpage{1780}
\lastpage{1826}

\makeatletter
\newcommand{\widetildes}{\tilde}

\newcommand{\eqref}[1]{(\ref{#1})}
\newtheorem{Lemma}{Lemma}[section]
\newtheorem{Proposition}[Lemma]{Proposition}

\newtheorem{Theorem}[Lemma]{Theorem}

\newproclaim{Problem}[Lemma]{Problem}
\newproclaim{Example}[Lemma]{Example}
\newproclaim{OP}[Lemma]{Open Problem}

\newtheorem{Corollary}[Lemma]{Corollary}
\newproclaim{Definition}[Lemma]{Definition}
\newproclaim{Condition}[Lemma]{Condition}

\newproclaim{Remark}[Lemma]{Remark}
\renewcommand{\mathbf}{\bolds}
\renewcommand{\P}{\mathbb{P}}
\newcommand{\prob}{\mathbb{P}}
\newcommand{\Pv}{\mathbb{P}}

\newcommand{\E}{\mathbb{E}}
\newcommand{\Ev}{\mathbb{E}}

\newcommand{\GG}{\mathcal{G}}
\newcommand{\NN}{\mathcal{N}}
\newcommand{\TT}{\mathcal{T}}
\newcommand{\Rbold}{{\mathbb{R}}}
\newcommand{\Nbold}{{\mathbb{N}}}

\newcommand{\expec}{\mathbb{E}}
\newcommand{\bfd}{\mathbf{d}}

\newcommand{\cC}{\mathcal{C}}
\newcommand{\cE}{\mathcal{E}}
\newcommand{\cF}{\mathcal{F}}
\newcommand{\cG}{\mathcal{G}}
\newcommand{\cL}{\mathcal{L}}
\newcommand{\cP}{\mathcal{P}}
\newcommand{\cV}{\mathcal{V}}

\newcommand{\bR}{\mathbb{R}}
\newcommand{\vC}{\mathbf{C}}
\newcommand{\vW}{\mathbf{W}}

\newcommand{\R}{\mathbb{R}}
\newcommand{\N}{\mathbb{N}}
\newcommand{\dd}{\mathrm{d}}

\newcommand{\Poi}{\operatorname{Poi}}

\newcommand{\CN}{\mathcal{N}}
\newcommand{\CT}{\mathcal{T}}

\newcommand{\one}{\mathbf{1}}
\newcommand{\wih}{\widehat}
\newcommand{\La}{\Lambda}
\newcommand{\ve}{\varepsilon}

\makeatother

\begin{document}
\begin{frontmatter}

\title{Degree distribution of shortest path trees
and bias of network sampling algorithms}
\runtitle{Degree distribution of shortest path trees}

\begin{aug}
\author[A]{\fnms{Shankar}~\snm{Bhamidi}\corref{}\ead[label=e1]{bhamidi@email.unc.edu}\thanksref{T1}},
\author[B]{\fnms{Jesse}~\snm{Goodman}\ead[label=e2]{jgoodman@tx.technion.ac.il}\thanksref{T2}},
\author[C]{\fnms{Remco}~\snm{van~der~Hofstad}\ead[label=e3]{rhofstad@win.tue.nl}\thanksref{T3}}
\and
\author[C]{\fnms{J\'ulia}~\snm{Komj\'athy}\ead[label=e4]{j.komjathy@tue.nl}\thanksref{T3}}
\thankstext{T1}{Supported in part by NSF  Grants DMS 1105581, 1310002 and NSF-SES Grant 1357606.}
\thankstext{T2}{The work was conducted partly while at
Eurandom, Leiden University, and Technion, and was supported in part
through ERC Advanced Grant 267356 VARIS. Also supported in part by
ISF Grant 915/12.}
\thankstext{T3}{Supported in part by The Netherlands
Organization for Scientific Research (NWO).}
\runauthor{Bhamidi, Goodman, van~der~Hofstad and Komj\'athy}
\affiliation{University of North Carolina, Technion, Eindhoven
University of Technology\break
and Eindhoven University of Technology}
\address[A]{S. Bhamidi\\
Department of Statistics\\
\quad and Operations Research\\
University of North Carolina\\
304 Hanes Hall\\
Chapel Hill, North Carolina 27516\\
USA\\
\printead{e1}} 
\address[B]{J. Goodman\\
Department of Mathematics\\
Technion\\
Haifa 32000\\
Israel\\
\printead{e2}}
\address[C]{R. van der Hofstad\\
J\'ulia Komj\'athy\\
Department of Mathematics\\
\quad and Computer Science\\
Eindhoven University of Technology\\
P.O. Box 513\\
5600 MB Eindhoven\\
The Netherlands\\
\printead{e3}\\
\phantom{E-mail:\ }\printead*{e4}}
\end{aug}

\received{\smonth{10} \syear{2013}}
\revised{\smonth{4} \syear{2014}}

%
\begin{abstract}
In this article, we explicitly derive the limiting degree distribution
of the shortest path tree from a single source on various random
network models with edge weights. We determine the {asymptotics of the
degree distribution for large degrees} of this tree and compare it
to the degree distribution of the original graph. We perform this
analysis for the complete graph with edge weights that are powers
of exponential random variables (weak disorder in the stochastic
mean-field model of distance), as well as on the configuration model
with edge-weights drawn according to any continuous distribution. In
the latter, the focus is on settings where the degrees obey a power
law, and we show that the shortest path tree again obeys a power law
with the \emph{same} degree power-law exponent.
We also consider random $r$-regular graphs for large $r$, and show that
the degree distribution of the shortest path tree is closely
related to the shortest path tree for the stochastic {mean-field}
model of distance. We use our results to {shed light on} an
empirically observed bias in network sampling methods.

This is part of a general program initiated in previous works by
Bhamidi, van der Hofstad and Hooghiemstra
[\emph{Ann. Appl. Probab.} \textbf{20} (2010) 1907--1965],
[\emph{Combin. Probab. Comput.} \textbf{20} (2011) 683--707],
[\emph{Adv. in Appl. Probab.} \textbf{42} (2010) 706--738]
of analyzing the effect of attaching
random edge lengths on the geometry of random network models.
\end{abstract}

%
\end{frontmatter}

\section{Introduction}\label{sint}

In the last few years, there has been an enormous\break amount of empirical
work in understanding properties of real-world networks, especially
data transmission networks such as the Internet. One functional which
has witnessed intense study and motivated an enormous amount of
literature is the \emph{degree distribution} of the network. Many
real-world networks are observed to have a \emph{heavy-tailed} degree
distribution. More precisely, empirical data suggest that if we look at
the empirical proportion ${\widehat p}_k$ of nodes with degree $k$, then
%
\begin{equation}
\label{eqEmpiricalDegreeTail} {\widehat p}_k \approx1/k^\tau,\qquad k\to\infty.
\end{equation}
The quantity $\tau$ is called the \textit{degree exponent} of the network
and plays an important role in predicting a wide variety of properties,
ranging from the typical distance between different nodes, robustness
and fragility of the network, to diffusion properties of viruses and
epidemics; see \cite
{Hofs10c,Durr06,ChuLu06c,newman,Dorogovtsev2003evolution,Newman-Barb-Watts}
and the references therein.

In practice, such network properties often cannot be directly measured
and are
estimated via indirect observations. The degree of a given node, or
whether two given
nodes are linked by an edge, may not be directly observable. One method
to overcome
this issue is to send probes from a single source node to every other
node in the
network, tracking the paths that these probes follow. This procedure,
known as
\emph{multicast}, gives partial information about the underlying
network, from which
the true structure of the network must be inferred; see, for example,
\cite{clauset-newman08,BroCla2001,faloutsos1999power,govindan2000heuristics,lakhina2003,pansiot1998routes}.

Probes sent between nodes to explore the structure of such networks are
assumed to follow \emph{shortest paths} in the following sense. These
networks are described not only by their graph structure but also by
costs or weights across edges, representing congestion across the edge
or economic costs for using {it}. The total weight of any given
path is the sum of edge weights along the path. Given a source node and
a destination node, a shortest path is a (potentially nonunique) path
joining these nodes with smallest total weight. It is generally
believed that the path that data actually takes is not the shortest path,
but that the shortest path is an acceptable approximation of the actual
path. For our models, the shortest paths between vertices will always
be unique.

For a given source node, the union of the shortest paths to all other
nodes of the network defines a subgraph of the underlying network,
representing the part of the network that can be inferred from the
multicast procedure. When all shortest paths are unique, which we
assume henceforth, this subgraph is a tree, called the \emph{shortest
path tree}. This will be the main object of study in this paper.

Given the shortest path tree and its degree distribution, one can then
attempt to infer the degree distribution of the whole network.
Empirical studies such as \cite{clauset-newman08,lakhina2003} show
that this may create a \emph{bias}, in the sense that the observed
degree distribution of the tree might differ significantly from the
degree distribution of the underlying network. Thus a theoretical
understanding of the shortest path tree, including its degree
distribution and the lengths of paths between typical nodes, is of
paramount interest.

By definition, the unique path in the shortest path tree from the
source $v_s$ to any given target vertex $v_t$ is the shortest path in
the weighted network between $v_s$ and $v_t$.
Thus the shortest path tree minimizes path lengths, not the total
weight of a spanning tree. Hence it is different from the \emph{minimal
spanning tree}, the tree for which the total weight over all edges is
the tree is minimal. The last few years have seen a lot of interest in
the statistical physics community for the study of disordered random
systems {that} bridge these two regimes, with models proposed to
interpolate between the shortest weight regime (first passage
percolation or weak disorder) and the minimal spanning tree regime
(strong disorder); see \cite{weak-strong-diso}. Consider a connected
graph $\cG_n = (\cV_n, \cE_n)$ on $n$ vertices with edge lengths
$\cL
_n:=\{l_e\dvtx e\in\cE_n\}$. Now fix {the} disorder parameter
$s\in
\bR_+$, change the edge weights to $\cL_n(s):=\{l_e^s\colon e\in
\cG_n\}$ and consider the shortest paths corresponding to the weights
$\cL
_n(s)$. For finite $s$, this is called the \emph{weak disorder}
regime. As $s\to\infty$, it is easy to check that the optimal path
between any two vertices converges to the path between these two
vertices in the minimal spanning tree where one uses the original edge
weights $\cL(n)$ to construct the minimal spanning tree. This is called
the \emph{strong disorder regime}. The parameter $s$ allows one to
interpolate between these two regimes. Understanding properties of the
shortest path tree and its dependence on the parameter $s$ is then of relevance.

The aim of this paper is to study the degree distribution of shortest
path trees, motivated by these questions from network sampling and
statistical physics.

\subsection{Mathematical model}
\label{secmodel-def}
In order to gain insight into these properties, we need to model (a)
the underlying networks and (b) the edge weights. We shall study two
main settings in this paper, the first motivated by network sampling
issues and the second to understand weak disorder models.

\begin{longlist}[(a)]
\item[(a)]\emph{Configuration model with arbitrary edge weights}: An
array of models have been proposed to capture the structure of
empirical networks, including preferential attachment-type models \cite
{AB99,BR04,BRST01} and, what is relevant to this study, the
configuration model ${\mathrm{CM}}_n(\mathbf{d})$ \cite
{Boll01,MolRee95} on $n$ vertices
given a degree sequence ${\mathbf{d}}=\mathbf
{d}_n=(d_1,\ldots,d_n)$ which is constructed as follows. Let
$[n]:=\{1,2,\ldots, n\}$ denote the vertex set of the graph. To vertex
$i\in
[n]$, attach $d_i$ half-edges, and write $\ell_n=\sum_{i\in[n]} d_i$
for the total degree, assumed to be even. (For $d_i$ drawn
independently from a common degree distribution $D$, $\ell_n$ may be
odd; if so, select one of the $d_i$ uniformly at random and increase it
by 1.) Number the half-edges in any arbitrary order from $1$ to $\ell
_n$, and sequentially pair them uniformly at random to form complete
edges. More precisely, at each stage pick an arbitrary unpaired
half-edge and pair it to another uniformly chosen unpaired half-edge to
form an edge. Once paired, remove the two half-edges from the set of
unpaired half-edges and continue the procedure until all half-edges are
paired. Call the resulting multi-graph ${\mathrm{CM}}_n(\mathbf{d})$.\looseness=-1

\begin{figure}[b]

\includegraphics{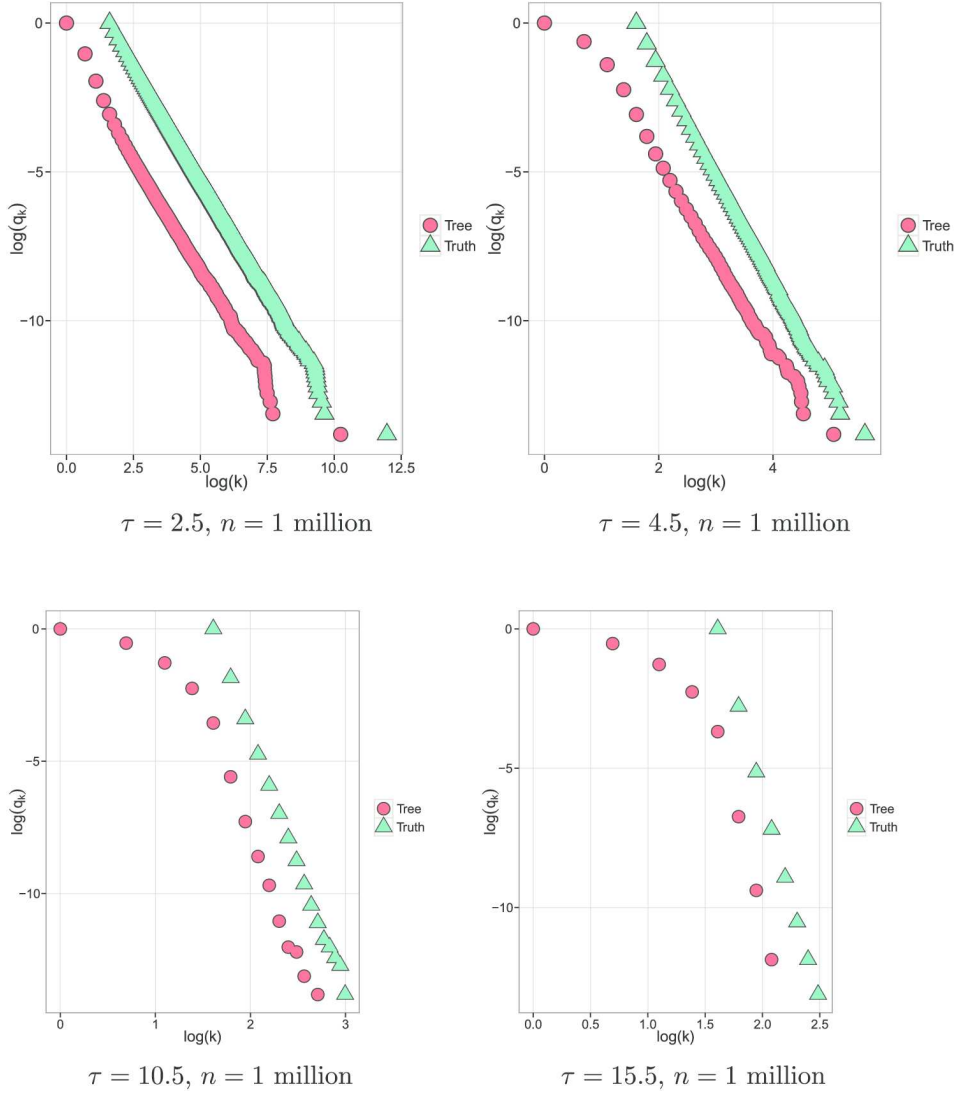}

\caption{Empirical distributions of underlying degrees (``truth'') in
the full graph and observed degrees (``tree'') in the shortest path
tree, shown in log--log scale. The vertical axis measures the tail
proportion $q_k=\sum_{j\geq k} \widehat{p}_j^{(n)}$ of vertices
having degree at least $k$. The underlying graphs are realizations of
the configuration model on $n$ vertices with power-law degree
distributions having exponent $\tau$ (and minimal degree $5$ so as to
ensure connectivity). Edge weights are i.i.d. exponential variables.}
\label{figsimulations}
\end{figure}

Although self-loops and multiple edges may occur, under mild conditions
on the degree sequence ${\mathbf{d}}$, these become rare as
$n\to\infty$; see, for example, \cite{Jans09} or \cite{Boll01} for
more precise results in this direction. For the edge weight
distribution, we will assume any continuous distribution with a
density. In the case of infinite-variance degrees, we need to make
stronger assumptions and only work with exponential edge weights and
independent and identically distributed (i.i.d.) degrees having a
power-law distribution.

\item[(b)]\emph{Weak disorder and the stochastic mean-field model}: The
complete graph can serve as an easy mean-field model for data
transmission, and for many observables it gives a reasonably good
approximation to the empirical data; see \cite{Mieghem09}. The complete
graph with random exponential mean one edge weights is often refered to
as the \emph{stochastic mean-field model of distance} and has been one
of the standard workhorses of probabilistic combinatorial optimization;
see \cite{janson123,aldous-assignment,aldous2003objective,wastlund} and
the references therein. In this context, we consider the weak disorder
model where, with $s > 0$ fixed, the edge lengths are i.i.d. copies of
$E^s$, where $E$ has an exponential distribution with mean one. In
\cite
{janson123}, the optimal paths were analyzed when $s=1$, and in \cite
{BhamidiHofstad2012}, the case of general $s$ was studied as a
mathematically tractable model of weak disorder.\looseness=-1
\end{longlist}

%
\subsection{Our contribution}
We rigorously analyze the asymptotic degree distribution of the
shortest path tree in the two settings described above. We give an
explicit probabilistic description of the limiting degree distribution
that is intimately connected to the random fluctuations of the length
of the optimal path. These in turn are intimately connected to
Bellman--Harris--Jagers continuous-time branching processes (CTBP)
describing local neighborhoods in these graphs. By analyzing these
random fluctuations, we prove that the limiting degree distribution has
markedly different behavior depending on the underlying graph:

\begin{longlist}[(ii)]
\item[(i)]\emph{Configuration model}:
The shortest path tree has the \emph{same degree exponent} $\tau$ as
the underlying graph for \emph{any} continuous edge weight distribution
when $\tau>3$, and for exponential mean one edge weights when $\tau
\in
(2,3)$. This reflects the fact that, for a vertex of unusually high
degree in the underlying graph, almost all of its adjoining edges (if
$\tau>3$) or a positive fraction of its adjoining edges (if $2<\tau<3$)
are likely to belong to the shortest path tree. See Figure~\ref{figsimulations}.

\begin{figure}

\includegraphics{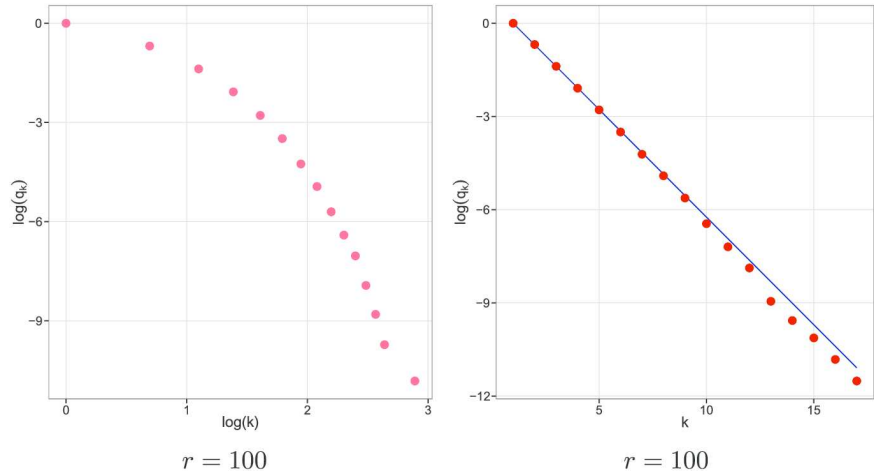}

\caption{Empirical distributions of observed degrees in the shortest path tree.
At left, both the degree $k$ and the tail proportion $q_k = \sum_{j
\geq k} \widehat p_j^{(n)}$ of vertices having degree at least
$k$ are shown in logarithmic scale; at right, only $q_k$ is shown in
logarithmic scale. The underlying graph is a random $r$-regular graph,
$r=100$, on $n=100\mbox{,}000$ vertices. The blue line in the right-hand graph
is the curve $q=2^{-k}$, corresponding to the $\operatorname{Geometric}(1/2)$
distribution. Edge weights are i.i.d. exponential variables.}
\label{figsimul-random-reg}
\end{figure}

\item[(ii)]\emph{Weak disorder}: Here the limiting degree distribution of
the shortest path tree has an exponential or stretched exponential tail
depending on the temperature $s$. Furthermore, this limiting degree
distribution arises as the limit $r\to\infty$ of the limiting degree
distribution for the $r$-regular graph when the edge weights are
exponential variables raised to the power $s$; see Figure~\ref{figsimul-random-reg} for the case $s=1$.
\end{longlist}

%
%
%
%
%
%
%
%
%
%

%
%
%


\subsection{Notation}
\label{secnot}
In stating our results, we shall write $v_s$ and $v_t$ for two vertices
(the ``source'' and the ``target'') chosen uniformly and independently
from a graph $\GG_n$ on vertex set $[n]=\{1,\ldots,n\}$, which will
either be the complete graph or a realization of the configuration
model. For the configuration model, we write $d_v$ for the degree of
vertex $v\in[n]$. On the edges of $\GG_n$, we place i.i.d. positive
edge weights $Y_e$ drawn from a continuous distribution. We denote by
$\TT_n$ the shortest path tree from vertex $v_s$, that is, the union
over all vertices $v\neq v_s$ of the (a.s. unique) optimal path from
$v_s$ to $v$. We write $\deg_{\TT_n}(v)$ for the degree of vertex $v$
in the shortest path tree and $\widehat{p}_k^{(n)}$ for the
proportion of vertices having degree $k$ in the shortest path tree.\looseness=1

We write $E$ for an exponential variable of mean $1$ and $\Lambda
\stackrel{d}{=}\log(1/E)$ for a standard Gumbel variable, that is,
$\P
(\Lambda
\leq x)=\exp(-{\mathrm e}^{-x})$.

\subsection{Organization of the paper}
We describe our results in Section~\ref{sResults} and set up the
necessary mathematical constructs for the proof in Section~\ref{sbackgroundmaterial}. Theorems about convergence of the degree
distribution have three parts:
\begin{itemize}[$\rhd$]
\item[$\rhd$] part (a) describes the \emph{limiting degree distribution
of a uniformly chosen vertex} in the shortest path tree; this is proved
in Section~\ref{spart-a-proofs};
\item[$\rhd$] part (b) states the convergence of the \emph{empirical
degree distribution} in the shortest path tree to the asserted limit
from part (a); this is proved in Section~\ref{spart-b-proofs};
\item[$\rhd$] part (c) identifies the limiting \emph{expected degree}
in the shortest path tree; this is proved in Section~\ref{spart-c-proofs}.
\end{itemize}
Section~\ref{sResults} also contains results about the tail behavior
of the degrees in the shortest path tree, proved in Section~\ref{ssasymptotics}, and a link between the limiting degree distributions
and those for breadth-first tree setting, proved in Section~\ref{sY=1Proof}.

\section{Main results and discussion}\label{sResults}
We now set out {to state} our main results.

\subsection{Weak disorder in the stochastic mean-field model}\label
{ssKnResults}

Let $\cG_n(s)$ denote the complete graph with each edge $e$ equipped
with an i.i.d. edge weight $l_e = E^s $ where $E\sim\exp(1)$ and $s>
0$. Here we describe our results for the shortest path tree $\mathcal
{T}_n:=\CT
_n(s)$ from a randomly selected vertex.
Let $E_i$, $i=1,2,\ldots,$ denote independent copies of $E$. Define
$X_1<X_2<\cdots$ by
%
\begin{equation}
\label{eqPPPEs} X_i=(E_1+\cdots+E_i)^s;
\end{equation}
equivalently, $(X_i)_{i\geq1}$ are the ordered points of a Poisson
point process with intensity measure
%
\begin{equation}
\label{eqmusDefinition} \mathrm{ d}\mu_s(x) = \frac{1}{s}x^{1/s-1}\,\mathrm{
d}x.
\end{equation}
Let $\Gamma(\cdot)$ be the Gamma function, and set
%
\begin{equation}
\label{eqlambdas} \lambda_s = \Gamma(1+1/s)^s;
\end{equation}
a short calculation verifies that $\int_0^\infty{\mathrm e}^{-\lambda
_s x}\,\dd\mu
_s(x)=1$. Then there exists a random variable $W$ with $W>0$ and $\E
(W)=1$ whose law is uniquely defined by the recursive distributional equation
%
\begin{equation}
\label{eqWfors} W \stackrel{d} {=}\sum_{i\geq1} {\mathrm
e}^{-\lambda_s X_i}W_i,
\end{equation}
where $W_1,W_2,\ldots$ are i.i.d. copies of $W$. This identity will
arise from the \emph{basic decomposition} of a certain continuous-time
branching process, and the uniqueness in law of $W$ follows from
standard arguments; see Section~\ref{sssBPapprox}.

Our first theorem describes the degrees in the shortest path tree for
the weak-disorder regime from Section~\ref{sint}:

\begin{Theorem}\label{tKnDegrees}
Let $s>0$, and place i.i.d. positive edge weights with distribution
$E^s$ on the edges of the complete graph $K_n$. Let $(X_i)_{i\geq1}$
be as in \eqref{eqPPPEs}, let $(\Lambda_i)_{i\geq1}$ be i.i.d. standard Gumbel variables and let $(W_i)_{i\geq1}$ be an i.i.d. sequence of copies of $W$. Then:
\begin{longlist}[(a)]
\item[(a)] the degree $\deg_{\TT_n}(V_n)$ of a
uniformly
chosen vertex in the shortest path tree converges in distribution to
the random variable $\widehat{D}$ defined by
%
\begin{eqnarray}
\label{eqKndegreecharacterisation} \widehat{D} = 1+\sum_{i\geq1}
\mathbh{1}_{\{\Lambda_i +\log W_i +
\lambda_s X_i < M\}}
\nonumber
\\[-8pt]
\\[-8pt]
\eqntext{\mbox{with}\quad \displaystyle M=\max_{i\in\N} (
\Lambda_i+\log W_i-\lambda_s X_i);}
\end{eqnarray}
\item[(b)] the empirical degree distribution in
the shortest path tree converges in probability as $n\rightarrow\infty$,
%
\begin{equation}
\label{eqKnempiricaldegree} \widehat p_k^{(n)}= \frac{1}n \sum
_{v\in[n]} \one_{\{\deg_{\TT
_n}(v)=k\}} \stackrel{{\mathbb P}} {
\longrightarrow}\Pv(\widehat D=k);
\end{equation}
\item[(c)]
the expected limiting degree is $2$, that is, as $n\rightarrow\infty$,
%
\begin{equation}
\label{eqKnaveragedegree} \Ev\bigl[\deg_{\TT_n}(V_n)\bigr] \to\Ev[
\widehat D]=2.
\end{equation}
\end{longlist}
\end{Theorem}

We remark that $\widehat{D}$ and $M$ take finite values: the law of
large numbers implies that $i^{-1} X_i\to1$ a.s., whereas $\Lambda
_i+\log W_i=o(i)$ a.s. by Markov's inequality and the Borel--Cantelli lemma.

Since $\TT_n$ is a tree on $n$ vertices, and $V_n$ is a uniformly
chosen vertex, $\E[\deg_{\TT_n}(V_n)]=2(1-1/n)\to2$ as $n\to\infty
$. The
convergence in \eqref{eqKnaveragedegree} is in this sense a
triviality. However, when combined with the convergence in distribution
of $\deg_{\TT_n}(V_n)$ to $\widehat{D}$ from part~(a),
the assertion of \eqref{eqKnaveragedegree} is that no mass is lost
in the limit; that is, the variables $\deg_{\TT_n}(V_n)$, $n\in\N$, are
uniformly integrable. In practical terms, this means that a small
number of vertices cannot carry a positive proportion of the degrees in
the shortest path tree.

The following theorem describes the tail of the degree distribution in
the tree in terms of the exponent $s$ on the exponential weights:

\begin{Theorem}\label{tKnDegreeTails}
Let $s>0$, and place i.i.d. positive edge weights with distribution
$E^s$ on the edges of the complete graph $K_n$:
\begin{longlist}[(a)]
\item[(a)]
For $s=1$, the variable $\widehat{D}$ defined by \eqref
{eqKndegreecharacterisation} is a geometric random variable with
parameter $\frac{1}{2}$. Then:
\item[(b)]
For $s<1$ and $k\to\infty$,
\[
\log\P ( \widehat{D}=k ) \sim-\lambda_s k^{s}.
\]

\item[(c)] For $s>1$ and $k\to\infty$,
\[
\log\P ( \widehat{D}=k ) \sim-(1-1/s)k\log k.
\]
\end{longlist}
\end{Theorem}

Theorem~\ref{tKnDegreeTails} shows that the tail asymptotics of
$\widehat{D}$ decay less rapidly when $s$ becomes small. Note that the
boundary case $s=0$ corresponds to \emph{constant} edge weights.
However, $\lambda_s\to\infty$ as $s{\searrow} 0$, and Theorem~\ref
{tKnDegreeTails} is not uniform over $s$. Indeed, the limit
$s{\searrow} 0$ is surprisingly subtle; see \cite{EGHN2013}.

\subsection{The configuration model with finite-variance
degrees}\label
{ssCMfiniteVarResults}

We next consider the configuration model for rather general degree
sequences $\bfd_n$, which may be either deterministic or random,
subject to the following convergence and integrability conditions. To
formulate these, we think of $\bfd_n$ as fixed, and choose a vertex
$V_n$ uniformly from $[n]$. We write $d_v$ for the degree of $v$ in the
original graph. Then the distribution of $d_{V_n}$ is the distribution
of the degree of a uniformly chosen vertex~$V_n$, conditional on the
degree sequence $\bfd_n$. We assume throughout that $d_v\geq2$ for
each $v\in[n]$.

%
\begin{Condition}[(Degree regularity)]\label{condCMFinVar}\label{cCMFinVar}
The degrees $d_{V_n}$ satisfy $d_{V_n}\geq2$ a.s. and, for some
random variable $D$ with $\prob(D>2)>0$ and $\E(D^2)<\infty$,
%
\begin{equation}
d_{V_n}\stackrel{d} {\longrightarrow}D, \qquad \E
\bigl(d_{V_n}^2\bigr) \to\E\bigl(D^2\bigr).
\end{equation}
Furthermore, the sequence $d_{V_n}^2 \log(d_{V_n})$ is uniformly
integrable. That is, for any sequence $a_n\to\infty$,
%
\begin{equation}
\limsup_{n\to\infty} \E\bigl(d_{V_n}^2
\log^+(d_{V_n}/a_n)\bigr) = 0.
\end{equation}
\end{Condition}

In the case where $\bfd_n$ is itself random, we require that the
convergences in Condition \ref{condCMFinVar} hold in probability. In
particular, Condition \ref{condCMFinVar} is satisfied when
$d_1,\ldots,d_n$ are i.i.d. copies of $D$ and $\E(D^2\log D)<\infty$.

{Define the distribution of the random variable $D^\star$---the \emph
{size-biased} version of~$D$---by}
%
\begin{equation}
\label{eqBiasedD} \prob\bigl(D^\star=k\bigr)=\frac{k\prob(D=k)}{\E(D)}.
\end{equation}
We define $\nu=\Ev(D^\star- 1)$; it is easily checked that $\nu=\E
[D(D-1)]/\E[D]$. The assumptions $d_{V_n}\geq2$ and $\prob(D>2)>0$
imply that $\nu>1$.

We take the edge weights to be i.i.d. copies of a random variable
$Y>0$ with a continuous distribution. Since $\nu>1$, we may define the
\emph{Malthusian parameter} $\lambda\in(0,\infty)$ by the
requirement that
%
\begin{equation}
\label{eqMalthusian} \nu\E\bigl({\mathrm e}^{-\lambda Y}\bigr)=1.
\end{equation}
Then there is a random variable $W$ whose law is uniquely defined by
the requirements that $W>0$, $\E(W)=1$ and
%
\begin{equation}
\label{eqWrecursionCMFinVar} W \stackrel{d} {=}\sum_{i=1}^{D^\star-1}
{\mathrm e}^{-\lambda Y_i} W_i,
\end{equation}
where $W_1,W_2,\ldots$ are i.i.d. copies of $W$.
Again, this identity is derived from the basic decomposition of a
certain branching process; see Section~\ref{sssBPapprox}.

The next theorem, the counterpart of Theorem~\ref{tKnDegrees}, is
about the degrees in the shortest path tree in the configuration model:

\begin{Theorem}\label{tCMtheoremfinitevar}
On the edges of the configuration model where the degree sequences
$\bfd
_n$ satisfy Condition \ref{condCMFinVar} with limiting degree
distribution $D$, place as edge weights i.i.d. copies of a random
variable $Y>0$ with a continuous distribution. Let $(\Lambda_i)_{i\geq
1}$, $(W_i)_{i\geq1}$ and $(Y_i)_{i\geq1}$ be i.i.d. copies of
$\Lambda$, $W$ and $Y$, respectively. Then:
\begin{longlist}[(a)]
\item[(a)] the degree $\deg_{\TT
_n}(V_n)$ of a
uniformly chosen vertex in the shortest path tree converges in
distribution to the random variable $\widehat{D}$ defined by
%
\begin{eqnarray}
\label{eqhatD-CMFinVarDefn} \widehat{D} = 1 + \sum_{i=1}^D
\mathbh{1}_{\{\Lambda_i + \log W_i +
\lambda Y_i < M\}}
\nonumber
\\[-8pt]
\\[-8pt]
\eqntext{\mbox{with }\displaystyle  M=\max_{1\leq i\leq D} (
\Lambda_i + \log W_i -\lambda Y_i );}
\end{eqnarray}
\item[(b)] the empirical degree
distribution in the shortest path tree converges in probability,
%
\begin{equation}
\label{eqCMtheoremfinitevarempiricaldegree} \widehat p_k^{(n)}= \frac{1}n \sum
_{v\in[n]} \one_{\{\deg_{\TT
_n}(v)=k\}} \stackrel{{\mathbb P}} {
\longrightarrow}\Pv(\widehat D=k)\qquad \mbox{as } n\to\infty;
\end{equation}
\item[(c)]
the expected limiting degree is $2$, that is, as $n\rightarrow\infty$,
%
\begin{equation}
\label{eqCMFinVaraveragedegree} \Ev\bigl[\deg_{\TT_n}(V_n)\bigr] \to\Ev[
\widehat D]=2.
\end{equation}
\end{longlist}
\end{Theorem}

As in Theorem~\ref{tKnDegrees}(c), part (c) implies that the degrees $\deg_{\TT
_n}(V_n)$, $n\in
\N
$, are uniformly integrable. Since $\deg_{\TT_n}(V_n)\leq d_{V_n}$,
part (c) follows from Condition~\ref{condCMFinVar}
using dominated convergence, but for completeness we will give a proof
that uses the stochastic representation \eqref{eqhatD-CMFinVarDefn} directly.

In \eqref{eqhatD-CMFinVarDefn}, the behavior of $\widehat{D}$ depends
strongly on the value of $D$, and in particular $\widehat{D}\leq D$
a.s. [This bound is clear in the original degree problem; to see it
from \eqref{eqhatD-CMFinVarDefn}, note that the summand for which
$M=\Lambda_i+\log W_i-\lambda Y_i$ must vanish.] Thus very large
observed degrees $\widehat{D}$ must arise from even larger original
degrees $D$. To understand this relationship, we define a family of
random variables $(\widehat{D}_k)_{k=1}^\infty$ by
%
\begin{eqnarray}
\label{eqhatDk-CMFinVarDefn} \widehat{D}_k = 1 + \sum_{i=1}^k
\mathbh{1}_{\{\Lambda_i + \log W_i
+ \lambda Y_i < M_k\}}
\nonumber
\\[-8pt]
\\[-8pt]
\eqntext{\mbox{with}\quad M_k=\displaystyle\max_{1\leq i\leq k}
( \Lambda_i + \log W_i -\lambda Y_i ).}
\end{eqnarray}
The distribution of $\widehat{D}_k$ corresponds to the limiting
distribution of $\deg_{\TT_n}(V_n)$ when, instead of being selected uniformly,
$V_n$ is conditioned to have degree $k$. The limiting distribution
$\widehat{D}$ from \eqref{eqhatD-CMFinVarDefn} is then the composition
%
\begin{equation}
\label{eqDhatcomposition} \widehat{D} \stackrel{d} {=} \widehat D_D,
\end{equation}
where $D$ has the asymptotic degree distribution from Condition \ref
{condCMFinVar}.

As well as depending on $k$, the distribution of $\widehat{D}_k$
depends on $\lambda>0$ and on the distributions of $(\Lambda
_i)_{i\geq
1}$, $(W_i)_{i\geq1}$ and $(Y_i)_{i\geq1}$, which we always assume to
be i.i.d. copies of $\Lambda$, $W$ and $Y$, respectively. We omit this
dependence from the notation.

The asymptotic behavior of $\widehat{D}$ is established by large values
of $D$, hence we study $\widehat{D}_k$ in the limit $k\to\infty$. The
following theorem shows that the form of \eqref{eqhatD-CMFinVarDefn}
and \eqref{eqhatDk-CMFinVarDefn} determines the asymptotic behavior
under very general conditions.

\begin{Theorem}\label{tCMdegreekfin}
Define $\widehat{D}_k$ according to \eqref{eqhatDk-CMFinVarDefn},
where the variables $(\Lambda_i)_{i\geq1}$, $(W_i)_{i\geq1}$ and
$(Y_i)_{i\geq1}$ are i.i.d. copies of arbitrary random variables
$\Lambda,W,Y$, independently for each $i\in\Nbold$, with $Y>0$ a.s. If
$\P(\Lambda>x)>0$ for each $x\in\R$, or if $\P(W>x)>0$ for each
$x\in\R
$, then $\widehat{D}_k=k(1-o_{\prob}(1))$ as $k\to\infty$.
\end{Theorem}

Theorem~\ref{tCMdegreekfin} shows that the proportion of summands in
\eqref{eqhatDk-CMFinVarDefn} that do not contribute to $\widehat{D}_k$
tends to $0$. In words, if the vertex has large degree in the original
graph, then it is likely that almost all of the outgoing edges will be
revealed by the shortest path tree.

On the contrary, the next result shows that under certain circumstances
the order of magnitude of the error is not necessarily small, that is,
finite behavior might modify the empirical data significantly compared
to the true limit behavior. We pay particular attention to the case
when the edge weights $(Y_i)_{i\geq1}$ are i.i.d. exponential or
uniform variables. In these cases we can determine the precise
asymptotic order of magnitude of the difference between the degrees in
the original graph and in the shortest path tree.

%
\begin{Theorem}\label{tRateOfConv}
Define $\widehat{D}_k,M_k$ according to \eqref{eqhatDk-CMFinVarDefn},
where the variables $(\Lambda_i)_{i\geq1}$, $(W_i)_{i\geq1}$ and
$(Y_i)_{i\geq1}$ are i.i.d. copies of a Gumbel variable $\Lambda$, a
positive random variable $W$ with $\E(W)<\infty$ and a positive random
variable $Y$.
Then:
\begin{longlist}[(a)]
\item[(a)]
$M_k=\log k + O_{\prob}(1)$ as $k\to\infty$;

\item[(b)]
if $\E({\mathrm e}^{\lambda Y})<\infty$, then $k-\widehat{D}_k$ is tight;

\item[(c)]
if $Y$ is a standard exponential variable and the Malthusian parameter
$\lambda$ satisfies $\lambda>1$, then
%
\begin{equation}
\label{eqleftedgesCM} k-\widehat{D}_k=\Theta_\prob
\bigl(k^{1-1/\lambda}\bigr);
\end{equation}

\item[(d)]
if $Y$ is a standard exponential variable and $\lambda=1$, then
%
\begin{equation}
k-\widehat{D}_k=\Theta_\prob(\log k).
\end{equation}
\end{longlist}
\end{Theorem}

Theorem~\ref{tRateOfConv}(b) applies to the
setting where $Y$ is a standard exponential variable and $0<\lambda<1$.
Interestingly, for the CM with exponential edge weights, one has
$\lambda=\nu-1$, where we recall that $\nu=\expec[D(D-1)]/\expec[D]$
denotes the expected forward degree. Thus $\lambda\in(0,1)$ precisely
when $\nu\in(1,2)$. The other cases are treated in
Theorem~\ref{tRateOfConv}(c) and~(d), where the behavior is really
different. Further, Theorem~\ref{tRateOfConv}(b) applies to the setting where $Y$ is a uniform
random variable, regardless of the value of $\lambda$.

An immediate consequence is the following corollary, handling the case
of i.i.d. degrees with power-law exponent $\tau>3$. Here we shall
assume that the distribution function $F(x)=\prob(D\leq x)$ of the
underlying degrees satisfies
%
\begin{equation}
\label{eqFcond} 1-F(x)=x^{1-\tau} L(x),
\end{equation}
where $x\mapsto L(x)$ is a slowly varying function as $x\to\infty$.

%
\begin{Corollary}\label{cD-hatDAsymp-CMFinVar}
Suppose that the configuration model degrees are i.i.d. copies of a
random variable $D$ whose distribution function satisfies \eqref
{eqFcond} with $\tau>3$. Then:
\begin{longlist}[(a)]
\item[(a)]
conditional on $\{D=k\}$, we have $\widehat{D}=D(1-o_{\prob
}(1))$ in the
limit $k\to\infty$, and
\item[(b)]
the distribution function of $\widehat{D}$ satisfies \eqref{eqFcond}
also, for the same $\tau$.
\end{longlist}
If in addition $\nu>2$ and the edge weights are exponentially
distributed, then
\begin{longlist}[(c)]
\item[(c)]
conditional on $\{D=k\}$, we have $D-\widehat{D} = \Theta_\prob
(k^{1-1/(\nu-1)})$ in the limit \mbox{$k\to\infty$}.
\end{longlist}
\end{Corollary}

Corollary~\ref{cD-hatDAsymp-CMFinVar}(a) and (b) show that large degrees are asymptotically fully detected
in the shortest path tree. Corollary~\ref{cD-hatDAsymp-CMFinVar}(c)
provides a {counterpart} by showing that $\widehat{D}$, though
asymptotically of the same order as $D$, may nevertheless be
substantially smaller when $D$ is of moderate size. Furthermore, this
effect is accentuated when $\nu$ is large.

Note that Theorems \ref{tCMdegreekfin}--\ref{tRateOfConv} and thus
Corollary~\ref{cD-hatDAsymp-CMFinVar} rely heavily on the fact that
the underlying degree distribution and the Malthusian parameter
$\lambda$
stay fixed whereas $k$ is large. In other words, these results pertain
to a single vertex of unusually large degree. In particular, Theorems
\ref{tCMdegreekfin}--\ref{tRateOfConv} do \emph{not} hold for the
random $k$-regular graph in the limit $k\to\infty$. In that case every
vertex---not just the target vertex---has degree $k$ and hence the
Malthusian parameter $\lambda=k-1$ tends to infinity \emph{together} with
the degree $k$. In the context of an $r$-regular graph, Theorems~\ref
{tCMdegreekfin}--\ref{tRateOfConv} apply instead to the asymptotic
degree behavior of a vertex of degree $k$ added artificially to the
random $r$-regular graph on $n$ vertices, with $r$ fixed, $k \gg r$ and
$n\to\infty$.

\subsection{The configuration model with infinite-variance
degrees}\label{ssCMinfiniteVar}

Section~\ref{ssCMfiniteVarResults} treats the configuration model with
degree distribution having a finite limiting variance. However, in many
real-life networks, this is not the case. Quite often, the available
empirical work suggests that the degrees in the network follow a
power-law distribution
with exponent $\tau\in(2,3)$.

Thus throughout this section we shall have in mind that the degrees
$d_1,\ldots,d_n$ of the
configuration model are i.i.d. copies of $D$, where $D\geq2$ a.s. and the distribution function $F(x)=\prob(D\leq x)$ satisfies \eqref{eqFcond}
for $2<\tau<3$ and $x\mapsto L(x)$ a slowly varying function as $x\to
\infty$. We further assume that the edge weights are standard
exponential random variables.

In the parameter range $2<\tau<3$, the degree distribution has finite
mean but infinite variance. Hence the size-biased distribution in
\eqref
{eqBiasedD} is well defined, but has infinite mean, and the Malthusian
parameter in \eqref{eqMalthusian} does not exist. Instead, let $V$ be
the positive ({nontrivial}) random variable
that satisfies
%
\begin{equation}
\label{eqexplosiontime-representation1} V \stackrel{d} {=} \sum_{i=1}^{\infty}
\frac{E_i}{1+\sum_{j=1}^i (D^{{\star}}_j-2)},
\end{equation}
where $E_i$ is an i.i.d. collection of exponential random variables,
and independently, $D_j^{\star}$ are i.i.d. copies of the
{size-biased} distribution defined in \eqref{eqBiasedD}, now having
infinite mean. It is not hard to see that $V$ also satisfies
%
\begin{equation}
\label{eqexplosiontime-representation} V \stackrel{d} {=} \min_{i=1, \ldots, D^\star-1} (
E_i+V_i ),
\end{equation}
where $E_i$ and $V_i$ are i.i.d. copies of $E$ and $V$, respectively.
This recursive characterization can be derived again from the basic
decomposition of Markov chains.

Our next theorem describes the behavior of degrees in the shortest path
tree on the configuration model with i.i.d. infinite-variance degrees
and exponential edge weights:

%
\begin{Theorem}\label{tCMdegreesinfvar}
On the edges of the configuration model whose degree sequence $\bfd_n$
is given by independent copies of $D$, where the distribution function
of $D$ satisfies \eqref{eqFcond} with $\tau\in(2,3)$, place i.i.d. edge weights distributed as $E$, an exponential random variable of mean
1. Let $(V_i)_{i\geq1}$ and $(E_i)_{i\geq1}$ be i.i.d. copies of $V$
and $E$, respectively. Then:
\begin{longlist}
\item[(a)]
the degree $\deg_{\TT_n}(V_n)$ of a uniformly chosen vertex in the shortest
path tree converges in distribution to the random variable $\widehat
{D}$ defined by
%
\begin{equation}
\label{eqhatD-CMInfVarDefn} \widehat{D}=1+\sum_{i=1}^D
\mathbh{1}_{\{V_i-E_i > \xi\}}\qquad\mbox{with } \xi= \min_{1\leq i\leq D}(V_i+E_i);
\end{equation}
\item[(b)] the empirical degree distribution
in the shortest path tree converges in probability, that is,
%
\begin{equation}
\label{eqCMinfvarempiricaldegree} \widehat p_k^{(n)}= \frac{1}n \sum
_{v\in[n]} \one_{\{\deg_{\TT
_n}(v)=k\}} \stackrel{{\mathbb P}} {
\longrightarrow}\Pv(\widehat D=k);
\end{equation}
\item[(c)]
the expected limiting degree is $2$, that is, as $n\rightarrow\infty$,
%
\begin{equation}
\label{eqCMInfVaraveragedegree} \Ev\bigl[\deg_{\TT_n}(V_n)\bigr] \to\Ev[
\widehat D]=2.
\end{equation}
\end{longlist}
\end{Theorem}

As with Theorem~\ref{tCMtheoremfinitevar}(c), {part (c) of
Theorem~\ref{tCMdegreesinfvar}} asserts that the degrees $\deg_{\TT_n}(V_n)
$, $n\in\N$, are uniformly integrable, and we will give both a
dominated convergence proof and a proof using the stochastic
representation \eqref{eqhatD-CMInfVarDefn}.

As in Section~\ref{ssCMfiniteVarResults}, we wish to understand the
asymptotic behavior of the degrees by looking at vertices with large
original degree. Thus, we define a family of random variables
$(\widehat
{D}_k)_{k=1}^\infty$ by
%
\begin{equation}
\label{eqhatDk-CMInfVarDefn} \widehat{D}_k=1+\sum_{i=1}^k
\mathbh{1}_{\{V_i-E_i > \xi_k\}} \qquad\mbox{with } \xi_k = \min
_{1\leq i\leq k}(V_i+E_i).
\end{equation}
Then the following theorem describes the degree in the shortest path
tree of a vertex conditioned to have a large original degree:

%
\begin{Theorem}\label{tCMdegreekinf}
Define $\widehat{D}_k$ according to \eqref{eqhatDk-CMInfVarDefn},
where the variables $(V_i)_{i\geq1}$ and $(E_i)_{i\geq1}$ independent
i.i.d. copies of arbitrary continuous positive random variables $V$
and $E$, respectively. If $\prob(V<x)$ and $\prob(E<x)$ are positive
for each $x>0$, then for $V,E$ independent, $p=\prob(V>E)$ satisfies
$0<p<1$ and,
as $k\to\infty$,
%
\begin{equation}
\label{eqbinasymp} \widehat{D}_k=p \cdot k \cdot\bigl(1+o_{\prob}(1)
\bigr).
\end{equation}
\end{Theorem}

Theorem~\ref{tCMdegreekinf} asserts that an asymptotic fraction $p$
(neither $0$ nor $1$) of the summands in \eqref{eqhatDk-CMInfVarDefn}
contribute to $\widehat{D}_k$. Compared to Theorem~\ref
{tCMdegreekfin}, where $p=1$, the difference stems from the fact
that $V$ is supported on $(0,\infty)$ whereas $\Lambda+\log W$ is
supported on $(-\infty,\infty)$.

%
\begin{Corollary}
If the distribution function of the configuration model degrees $D$
satisfies \eqref{eqFcond} with $\tau\in(2,3)$, then:
\begin{longlist}[(a)]
\item[(a)] conditional on $\{D=k\}$, we have $\widehat{D}=p \cdot D
\cdot(1+o_{\prob}(1))$ in the limit $k\to\infty$, and
\item[(b)] the distribution function of $\widehat{D}$ satisfies
\eqref
{eqFcond} also, for the same $\tau$.
\end{longlist}
\end{Corollary}
%

\subsection{Discussion}
In this section, we discuss our results and compare them to existing literature.

\subsubsection{Convergence to the limiting degree distribution}
Part (a) of Theorems~\ref{tKnDegrees}, \ref{tCMtheoremfinitevar} and
\ref{tCMdegreesinfvar} states that the degree distribution of a
single uniformly selected vertex converges to the distribution of
$\widehat{D}$. Part (b) strengthens this to state that the empirical
degree distribution converges in probability; that is, the (random)
proportion of vertices of degree $k$ in the shortest path tree $\TT_n$
is with high probability close to the limiting value $\prob(\widehat
{D}=k)$, for all $k$. Finally, part (c) states that the convergence of
the degree distribution from part (a) also happens in expectation.

Note that these convergences are not uniform over the choice of
original degree distribution or edge weight distribution; see the
remarks following Theorem~\ref{tKnDegreeTails} and Corollary~\ref
{cD-hatDAsymp-CMFinVar}.

\subsubsection{Degree exponents, bias and the effect of
randomness}\label{ssdegexp}

If the initial graph is the configuration model whose original degrees
obey a power law with exponent~$\tau$, then Theorems~\ref
{tCMdegreekfin} and \ref{tCMdegreekinf} show that
in both cases the power-law exponent $\tau$ is preserved via the
shortest path tree sampling procedure.\setcounter{footnote}{3}\footnote{To be precise, this is
proved only for $\tau>3$, for certain parts of the regime $\tau=3$ and
for $2<\tau<3$ with exponential edge weights.} In particular, if the
degrees from a shortest path tree are used to infer the power-law
exponent $\tau$, then asymptotically they will do so correctly.

In the literature, several papers consider the question of \emph{bias}.
Namely, do the observed degrees arising from network algorithms
accurately reflect the true underlying degree distribution, or can they
exhibit power law behavior with a modified or spurious exponent~$\widetilde
{\tau}$? This question has drawn particular attention in the setting of
the \emph{breadth-first search tree} (BFST), where paths are explored
in breadth-first order according to their number of edges, instead of
according to their total edge weight. Exact analysis \cite
{clauset-newman08} and numerical simulations \cite{lakhina2003} have
shown that the BFST can produce an apparent bias, in the sense that
observed degree distributions appear to follow a power law, for a
relatively wide range of degrees, when the true distribution does not.
Surprisingly, this phenomenon occurs even in the random $r$-regular
graph, where all vertices have degree $r$: defining
%
\begin{equation}
a_k^{(r)}=\frac{\Gamma(r)\Gamma(k-1+1/(r-2))}{(r-2)\Gamma
(r+1/(r-2))\Gamma(k)},
\end{equation}
the limiting degree distribution $\widehat{D}^{\mathrm{BFST}}$ satisfies
%
\begin{eqnarray}
\label{eqrBFSTtail} \prob\bigl(\widehat{D}^{\mathrm{BFST}}=k\bigr)&=&
a_k^{(r)}\qquad
\mbox{if }1\leq k\leq r\quad \mbox{and}
\nonumber
\\[-8pt]
\\[-8pt]
\nonumber
 a_k^{(r)}&\approx&
\frac{1}{r\cdot
k^{1-1/(r-2)}}\qquad\mbox{as }k\to\infty.
\end{eqnarray}
(See \cite{clauset-newman08}, Section~6.1; note that the requirement
$k\leq r$ is not mentioned in their discussion.) In this case, since
the underlying degrees are bounded, the power law in \eqref
{eqrBFSTtail} is of course truncated, and is therefore not a power law
in the sense of~\eqref{eqEmpiricalDegreeTail} or \eqref{eqFcond}.

The breadth-first search tree corresponds in our setup to the
nonrandom case where all edge weights are $1$. Although our proof of
Theorem~\ref{tCMtheoremfinitevar} relies on a continuous edge\vspace*{1pt} weight
distribution, we may nevertheless set $Y=1$ in the definition~\eqref
{eqhatD-CMFinVarDefn} of $\widehat{D}$. In this case, we recover the
limiting degree distribution arising from the breadth-first search tree:

%
\begin{Theorem}\label{tY=1}
Let $D$ be any degree distribution with $D\geq3$ a.s. and $\E
(D^2\log
D)<\infty$, and set $Y=1$. Then with $\lambda$ and $W$ as in Section~\ref{ssCMfiniteVarResults}, the limiting degree distribution
$\widehat
{D}$ from \eqref{eqhatD-CMFinVarDefn} is equal to the limiting degree
distribution for the breadth-first search tree identified in \cite{clauset-newman08},
Theorem~2.
\end{Theorem}

In particular, Theorem~\ref{tRateOfConv} (which makes no assumptions
on the edge weights except positivity) applies to the breadth-first
search tree degrees. Consequently, Theorem~\ref{tRateOfConv} and
Corollary~\ref{cD-hatDAsymp-CMFinVar} must be understood with the
caveat that they pertain to true power laws, but not truncated power
laws such as \eqref{eqrBFSTtail}.

For the truncated power law in \eqref{eqrBFSTtail} to look
convincingly like a true power law, $r$ must be relatively large. It is
worth noting, however, that the limiting degree distribution is \emph
{ill-behaved} in the limit $r\to\infty$: we have $a_1^{(r)}\to1$
and $a_k^{(r)}\to0$ for $k\geq2$, so that the degree of a typical
vertex converges to 1 and most vertices are leaves. In particular, the
truncated power law in \eqref{eqrBFSTtail} disappears in this limit.
Furthermore, the expected limiting degree $\E(\widehat{D}^{\mathrm
{BFST}})$ (which continues to be $2$ for each finite $r$) is reduced to
1 after taking $r\to\infty$, so that $\widehat{D}^{\mathrm
{BFST}}$ is
not uniformly integrable in this limit.

By way of comparison, the limiting degree distribution for the random
$r$-regular graph with i.i.d. exponential edge weights (perhaps raised
to some power $s>0$) is \emph{well behaved} in the limit $r\to\infty$,
and indeed converges\footnote{This follows from the convergence of the
collection ${(r^{-s}Y_i)_{i\in[r]}}$ of rescaled edge weights
toward the Poisson point process $(X_i)_{i=1}^\infty$ [cf. \eqref
{eqktuples-conv} and the surrounding material] and the consequent
convergence of the corresponding martingale limits $W$. Problems
related to the unbounded number of terms in \eqref
{eqKndegreecharacterisation} and \eqref{eqhatD-CMFinVarDefn} can be
handled by the observation that the collection ${(r^{-s}Y_i)_{i\in
[r]}}$ is stochastically dominated by ${(X_i)_{i\geq1}}$ for each
$r$.} to the limiting degree distributions for the complete graph
defined in Section~\ref{ssKnResults}. By Theorem~\ref
{tKnDegreeTails}, the tails of this distribution decay faster than a
power law, for any $s>0$.

Figure~\ref{figsimul-random-reg} shows a simulation of the case
$r=100$, $s=1$, with $n=10\mbox{,}000$. The observed degree distribution does
not resemble a power law at all, and in fact it agrees very closely
with the $\operatorname{Geometric}(1/2)$ distribution which, by Theorem~\ref
{tKnDegreeTails}(a) and the preceding discussion,
corresponds to the case $r\to\infty$. While not a proof, this strongly
suggests that the truncated power laws found in \cite
{clauset-newman08,lakhina2003} are anomalous and reflect specific
choices in the breadth-first search model. It would be of great
interest to understand under what conditions truncated power laws can
be expected to appear in general. It is tempting to conjecture that
spurious power laws do not arise whenever the edge weights are random
with support reaching all the way to $0$.

\subsubsection{Special cases}\label{sssSpecialCases}

The statement of Theorem~\ref{tKnDegreeTails} for $s=1$ is well known,
since in this case the shortest path tree is the \emph{uniform
recursive tree}, and the degrees in the uniform recursive tree can be
understood via martingale methods; see, for instance, \cite{Hofs10c}, Exercise
8.15, Theorem~8.2. The proof we give here is different, with
the main advantage that it is easier to generalize to the case $s\neq
1$. It is based on the representation \eqref
{eqKndegreecharacterisation} for $\widehat{D}$ together with the
observation that the martingale limit $W$ is a standard exponential
variable; see, for instance, \cite{Hofs10c} or \cite{janson123}, or
verify directly that $E$ satisfies \eqref{eqWfors}.

The $r$-regular graph on $n$ vertices corresponds to the choice $D=r$
in Theorem~\ref{tCMtheoremfinitevar}. If in addition the edge weights
are exponential, the martingale limit $W$ can be identified as a
Gamma$(\frac{r-1}{r-2},\frac{r-2}{r-1})$ random variable, that is, the
variable with Laplace transform $\phi_{W}(u)=(1+\frac
{r-2}{r-1}u)^{-(r-1)/(r-2)}$. Even though we can characterize $W$,
obtaining an explicit description of the law of $\widehat{D}$ (e.g.,
through its generating function) appears difficult.

\subsubsection{Branching processes: Limit random variables
\texorpdfstring{$W$ and $V$}{W and V}}\label{sssBPapprox}

In analyzing the shortest path tree $\TT_n$, it is natural to consider
the \emph{exploration process}, or first passage percolation, that
discovers $\TT_n$ gradually according to the distance from the source
vertex $v_s$. Starting from the subgraph consisting of $v_s$ alone,
reveal the original degree $d_{v_s}$. Reveal whether any of the
$d_{v_s}$ half-edges associated to $v_s$ form self-loops; if any do,
remove them from consideration. (This step is unnecessary in the
complete graph case.) For each remaining half-edge, there is an i.i.d. copy of the edge weight $Y$. Set $t_0=0$. Iteratively, having
constructed the subgraph with $i$ vertices and $i-1$ edges, wait until
the first time $t_i>t_{i-1}$ when some new vertex $v_i$ can be reached
from $v_s$ by a path of length $t_i$. (Thus $t_1$ will be equal to the
smallest edge weight incident to $v_s$, apart from self-loops.) Reveal
the degree $d_{v_i}$ and add the unique new edge in the path between
$v_i$ and $v_s$, using one of the $d_{v_i}$ half-edges associated to
$v_i$. For the remaining $d_{v_i}-1$ half-edges, remove any that form
self-loops or that connect to already explored vertices, and iterate
this procedure as long as possible. The subgraph so constructed will be
$\TT_n$.

When $n\to\infty$, no half-edge will form a self-loop or connect to a
previously explored vertex by any fixed stage $i$ of the exploration,
for any fixed $i$. It follows that the exploration process is well
approximated (at least initially) by a \emph{continuous-time branching
process} (CTBP) that we now describe.

Consider first the configuration model. The vertex $v_s$ is uniformly
chosen by assumption. The vertex $v_1$, however, is generally \emph
{not} uniformly chosen. Conditional on $v_s$ we have
%
\begin{equation}
\mathbb{P}(v_1=v \vert v_s) =
\frac{d_v \mathbh{1}_{\{v\neq v_s\}
}}{\sum_{w\neq
v_s} {d}_w}.
\end{equation}
(Note, e.g., that $d_{v_1}$ can never be $0$). Owing to the
finite mean assumption on the CM degrees, it follows that $\sum_{w\neq
v_s} \,d_w \sim n\E(D)$ and $\mathbb{P}(d_{v_1}=k \vert v_s)\approx
k\prob
(D=k)/\E
(D)$ in the limit $n\to\infty$. This \emph{size-biasing} effect means
that the number $d_{v_i}-1$ of new half-edges will asymptotically have
the distribution {$D^\star-1$, where $D^\star$ is} defined in
\eqref
{eqBiasedD}. The CTBP approximation for the CM is therefore the
following: An individual $v$ born at time $T_v$ has a random finite
number $N_v$ of offspring, born at times $T_v+Y_{v,1},\ldots,T_v+Y_{v,N_v}$. The $Y_{v,i}$ are i.i.d. copies of $Y$; the initial
individual $v_s$ has family size $N_{v_s} = d_{v_s}$; and all other
individuals have family size $N_v\stackrel{d}{=}D^\star-1$.

For the complete graph, the degrees are deterministic but large, and it
is necessary to rescale the edge weights: the collection of edge
weights incident to a vertex, multiplied by $n^s$, converges toward the
Poisson point process ${(X_i)_{i\geq1}}$ defined in~\eqref
{eqPPPEs}, for a formal version of this statement; see \eqref
{eqktuples-conv} below. The corresponding CTBP is as follows: Every
individual $v$ born at time $T_v$ has an infinite number of offspring,
born at times $T_v+X_{v,1},T_v+X_{v,2},\ldots,$ where
${(X_{v,i})_{i\geq1}}$ are i.i.d. copies of the Poisson point process
defined in \eqref{eqPPPEs}.

The random variables $W$ and $V$ from Sections~\ref{ssKnResults}--\ref
{ssCMinfiniteVar} arise naturally from these CTBPs. In the complete
graph context from Section~\ref{ssKnResults}--\ref
{ssCMfiniteVarResults}, the CTBPs grow exponentially in time, with
asymptotic population size $cW{\mathrm e}^{\lambda t}$ for $\lambda
=\lambda_s$
defined by~\eqref{eqlambdas} and $c>0$ a constant, and indeed $W$
arises as a suitable martingale limit; see~\cite{athreya}. For the CM
contexts from Sections~\ref{ssCMfiniteVarResults}--\ref
{ssCMinfiniteVar}, we must take the initial individual $v_s$ to have
degree distribution $D^\star-1$ in order to obtain the variables $W$
and $V$ (instead of $\widehat{W}$ and $\widehat{V}$ from Section~\ref{sbackgroundmaterial} below). When the family sizes $D^\star-1$ have
finite mean, as in Section~\ref{ssCMfiniteVarResults}, the population
size again grows asymptotically as $cW{\mathrm e}^{\lambda t}$ for
$\lambda$
given by \eqref{eqMalthusian}. In the setting of Section~\ref{ssCMinfiniteVar}, the CTBP explodes in finite time; that is, there is
an a.s. finite time $V=\lim_{k\to\infty} t_k$ at which the population
size diverges; see \cite{Grey73}. The recursive relations \eqref
{eqWfors}, \eqref{eqWrecursionCMFinVar} and \eqref
{eqexplosiontime-representation} result from conditioning on the size
and birth times of the first generation in the CTBP. For the uniqueness
in law of $W$, see \cite{levinson1960}, Theorem~4.1, page 111. 

We note that in all cases, the value of $W$ or $V$ is determined from
the initial growth of the branching process approximations: we can
obtain an arbitrarily accurate guess, with probability arbitrarily
close to $1$, by examining the CTBP until it reaches a sufficiently
large but finite size. In terms of the exploration process, it is
sufficient to examine a large but finite neighborhood of the initial
vertex. Large values of $W$ and small values of $V$ correspond to
faster than usual growth during this initial period, and thereafter the
growth is essentially deterministic.

In Theorems \ref{tKnDegrees} and \ref{tCMtheoremfinitevar}, a large
value of $M$ might be expected to correspond to one large value of
$W_i$, and a large value of $\widehat{D}$ might be expected to arise
from having many vertices $j$ with small values of $W_j$. As we shall
see in the proofs, however, this intuition is incorrect, and it is the
variables $\Lambda_i$, and secondarily the edge weights $Y_i$, whose
deviations are most relevant to the sizes of $M$ and $\widehat{D}$.

\subsubsection{Shortest path trees and giant components}

In Theorems \ref{tCMtheoremfinitevar} and \ref{tCMdegreesinfvar},
the hypothesis $D\geq2$ implies that $v_s$ and $v_t$ are connected
with high probability. If degrees $1$ or $0$ are possible, we must
impose the additional assumption that $\nu>1$ in Theorem~\ref
{tCMtheoremfinitevar}. Having made this assumption, the CM will have
a \emph{giant component}; that is, asymptotically, the largest
component will contain a fixed positive fraction of all vertices, and
the next largest component will contain $o(n)$ vertices. The variable
$W$ from Section~\ref{ssCMfiniteVarResults} has a positive probability
of being $0$, in which case we set $\log W=-\infty$, and the variable
$V$ from Section~\ref{ssCMinfiniteVar} has a positive probability of
being $\infty$. Furthermore, there is a positive probability that $\TT
_n$ contains only a fixed finite number of vertices, corresponding to
the case where the branching process approximations from Section~\ref{sssBPapprox} go extinct. (This possibility will be reflected
mathematically in the possibilities that $\widehat{W}_s=0$ in
Proposition~\ref{propWeightFinVarJointBiased} or $\widehat
{V}_s=\infty
$ in Proposition~\ref{propWeightInfVarJointBiased}.)

If we condition $v_s$ to lie in the giant component (corresponding to
nonextinction of the branching process started from $v_s$), then in
the resulting shortest path tree, the \emph{outdegree} of $v_t$ has the
same limiting conditional distribution as $\widehat{D}-1$ in Theorems
\ref{tCMtheoremfinitevar} and \ref{tCMdegreesinfvar}. The variable
$M$ (resp., $\xi$) equals $-\infty$ (resp., $\infty$)
whenever $W_i=0$ (resp., $V_i=\infty$) for each $i=1,\ldots,D$,
corresponding to the case that $v_t$ does not belong to the giant
component, and in this case the outdegree and the degree of $v_t$ are
both $0$.

\subsubsection{Open problems}
There are several interesting questions that serve as extensions of our results.
First, as discussed in Section~\ref{ssdegexp}, our results reveal the
existence or nonexistence of true power laws, but not truncated power
laws. A precise characterization of {when} truncated power laws
{arise} would be of great interest.

Second, many real-life networks have power{-}law behavior with
degree exponent $\tau\in(2,3)$. In this regime where the degrees have
infinite variance (as well as part of the regime $\tau=3$ when
Condition \ref{cCMFinVar} is not satisfied), it is natural to extend
beyond the exponential edge weights that we consider. We expect that
Theorems~\ref{tCMdegreesinfvar} and \ref{tCMdegreekinf} remain
valid with slight modification if the corresponding CTBP is explosive,
that is, if the CTBP reaches an infinite population in finite time.
When the corresponding\vspace*{1pt} CTBP is not explosive, even the probabilistic
form of the limiting distribution $\widehat{D}$ is unknown. Such a
representation would in particular be expected to give rise to the
limiting BFST degree distribution, as in Theorem~\ref{tY=1}.

Finally, real-world traceroute sampling typically uses more than just a
single source. It is natural to extend our model to several shortest
path trees from different sources. In this setup, the resulting
behavior might depend on whether we observe, for a given target vertex,
either the degree in each shortest path tree; or the degree in the
union of all shortest path trees; or the entire collection of incident
edges in each shortest path tree. In any of these formulations, we may
ask how accurately the observed degree reflects the true degree when
the number of sources is large, and whether this accuracy varies when
both the true degree and the number of sources are large.

\section{Limit theorems for shortest paths}\label{sbackgroundmaterial}
The proofs of Theorems~\ref{tKnDegrees}, \ref{tCMtheoremfinitevar}
and~\ref{tCMdegreesinfvar} are based on Propositions \ref
{propWeightKnJoint}, \ref{propWeightFinVarJointBiased} and \ref
{propWeightInfVarJointBiased}, respectively, which in turn follow from
\cite{BhamidiHofstad2012}, Theorem~1.1, \cite{BvdHHUnivFPP},
Theorems~1.2--1.3 and \cite{BhaHofHoo09b}, Theorem~3.2,
respectively.
These theorems determine the distribution of the shortest paths between
two uniformly chosen vertices in the complete graph, and in the
configuration model.
Since we need the results about shortest paths jointly across a
collection of several target points (i.e., not just between two
vertices), we state only the needed versions here. These modifications
easily follow from the results mentioned earlier, combined with an idea
about marginal convergence from the work of Salez \cite{Sale13} who
proved the joint convergence of typical distances between several
points for the particular case of the random $r$-regular graph with
exponential mean one edge lengths. His argument extends, however, to
the more general situation as well. We give an idea of how these
results were proven in Section~\ref{secidea-proof} but omit full proofs.
Our first proposition is about the joint convergence of shortest
weight paths on the complete graph. Recall the notation for $W$ from
Section~\ref{ssKnResults}.

%
\begin{Proposition}\label{propWeightKnJoint}
Consider the complete graph with edge weights distributed as $E^s$, $s>0$.
Let $v_1,\ldots,v_k$ be distinct vertices, all distinct from
$v_{\tilde
s}$ (the \emph{source vertex}), and denote the length of the shortest
path between $v_{\widetildes s}, v_i$ by $C_n(v_{\widetildes s},v_i)$. Then
%
\begin{equation}
 \bigl( \lambda_s n^s C_n(v_{\widetildes s},v_i)
- \log n \bigr)_{i=1}^k \stackrel{d} {\longrightarrow} ( -
\Lambda_i-\log W_{\widetildes
s} - \log W_i
)_{i=1}^k,
\end{equation}
where $\Lambda_1,\ldots,\Lambda_k$ are i.i.d. copies of $\Lambda$ and
$W_{\widetildes s},W_1,\ldots,W_k$ are i.i.d. copies of the random
variable $W$ from Section~\ref{ssKnResults}.
\end{Proposition}

Note that, due to the presence of the term $\log W_{\widetildes s}$, the
limiting variables in Proposition~\ref{propWeightKnJoint} are
exchangeable but not independent for different $i$. When $k=1$, the
case $s=1$ is due to \cite{janson123}, and the case $s\neq1$ is due to
\cite{BhamidiHofstad2012}.

For the configuration model with finite-variance degrees, we will need
to apply a similar result to the neighbors of the uniformly chosen
vertex $v_t$.

\begin{Proposition}\label{propWeightFinVarJointBiased}
Consider the configuration model with degrees satisfying Condition \ref
{condCMFinVar}.
Let $v_1,\ldots,v_k$ be distinct vertices, all distinct from
${v_{\widetildes s}}$, which may be randomly chosen but whose choice is
independent of the configuration model and of the edge weights. If the
degrees $(d_{v_1},\ldots,d_{v_k})$ converge jointly in distribution to
independent copies of $D^\star-1$, then there is a constant $\lambda>0$
and a sequence $\lambda_n\to\lambda$ such that
%
\begin{equation}
 \bigl( \lambda_n C_n({v_{\widetildes s}},v_i)
-\log n \bigr)_{i=1}^k \stackrel{d} {\longrightarrow} ( -
\Lambda_i-\log \widehat{W}_s - \log W_i + c
)_{i=1}^k,
\end{equation}
jointly in $i=1,\ldots,k$, where $c$ is a constant, $\Lambda_i$ are
i.i.d. copies of $\Lambda$, $W_1,\ldots,W_k$ are i.i.d. copies of the
variable $W$ from Section~\ref{ssCMfiniteVarResults} and $\widehat
{W}_s$ is a positive random variable, all independent of one another.
\end{Proposition}

{As discussed in Section~\ref{sssBPapprox},} each time we connect
a half edge of $v_t$ to another vertex, the probability of picking a
vertex of degree $k$ is approximately proportional to $k\cdot\Pv
(D=k)${. Thus, for each neighbor, the degree converges in distribution
to the size-biased variable $D^\star$ defined in \eqref{eqBiasedD},
and the number of half-edges not connected to $v_t$ converges in
distribution to $D^\star-1$. This motivates the assumptions on the
degrees in Proposition~\ref{propWeightFinVarJointBiased}.}

The constant $c$ arises as a function of the stable age-distribution of
the associated branching process \cite{BvdHHUnivFPP}. Since it does not
play a role in the proof, we omit a full description of this constant.

Finally, we state the corresponding result for the infinite-variance case:

\begin{Proposition}\label{propWeightInfVarJointBiased}
Consider the configuration model with i.i.d. degrees satisfying \eqref
{eqFcond} with $\tau\in(2,3)$.
Let $v_1,\ldots,v_k$ be distinct vertices, all distinct from
${v_{\widetildes s}}$, which may be randomly chosen but whose choice is
independent of the configuration model and of the edge weights. If the
degrees $(d_{v_1},\ldots,d_{v_k})$ converge jointly in distribution to
independent copies of the size-biased distribution $D^\star-1$, then
\[
\bigl( C_n({v_{\widetildes s}}, v_i)
\bigr)_{i=1}^k \stackrel {d} {\longrightarrow} (
\widehat{V}_s+V_i )_{i=1}^k,
\]
where $(V_i)_{i\geq1}$ are i.i.d. copies of the random variable $V$
from Section~\ref{ssCMinfiniteVar}, and $\wih V_s$ is a random
variable independent of $V_1,\ldots,V_k$.
\end{Proposition}

\subsection{Idea of the proof}
\label{secidea-proof}

We give the idea behind the proof of Proposition~\ref
{propWeightFinVarJointBiased}. The proofs of the other propositions
are similar, using the corresponding branching process approximations
of local neighborhoods as described in Section~\ref{sssBPapprox}.

Let ${(\mathbf{d}_n)_{n\geq1}}$ be a degree sequence
satisfying Condition \ref{condCMFinVar}, and fix a continuous positive
random variable $Y$. Let $\cG_n = ([n], \cE_n)$ be the configuration
model constructed from this degree sequence, with $\cE_n$ denoting the
edge set of the graph, and let the edge weights $\{Y_e\dvtx e\in\cE_n\}$
be i.i.d. copies of $Y$.

As in \eqref{eqBiasedD}--\eqref{eqMalthusian}, we define $\prob
(D_n^\star=k)=k\prob(d_{V_n}=k)/\E(d_{V_n})$ (the size-biasing of
$d_{V_n}$) and the corresponding size-biased expectations $\nu_n=\E
(D_n^\star-1)$, Malthusian parameters $\lambda_n$ satisfying $\nu
_n\E({\mathrm e}
^{-\lambda_n Y})=1$ and martingale limit $W^{(n)}$ satisfying
$W^{(n)} \stackrel{d}{=}\sum_{i=1}^{D_n^\star-1} {\mathrm
e}^{-\lambda
_n Y_i}
W_i^{(n)}$. Assuming Condition \ref{condCMFinVar}, we have $\nu
_n\to\nu$ (so that $\nu_n>1$ and $\lambda_n$, $W^{(n)}$ are
well defined for $n$ sufficiently large), $\lambda_n\to\lambda$ and
$W^{(n)}\stackrel{d}{\longrightarrow}W$.

\subsubsection{One target vertex: The case $k=1$}

Let us first summarize the ideas behind \cite{BvdHHUnivFPP},
Theorems~1.2--1.3, which derive the asymptotics for the
length of the optimal path between two selected vertices $v_0, v_1 \in
\cG_n$. To understand this optimal path, think of a fluid flowing at
rate one through the network using the edge lengths, started \emph
{simultaneously} from the two vertices $v_0, v_1$ at time $t=0$. When
the two flows collide, say at time $\Xi_n^{(1)}$, there exists one
vertex in both flow clusters. This implies that the optimal path is
created, and the length of the optimal path is essentially $2 \Xi
_n^{(1)}$.

Write ${(\cF_i(t))_{t\geq0}}$ for the flow process emanating from
vertex $v_i$. As described in Section~\ref{sssBPapprox}, these flow
processes can be approximated by independent Bellman--Harris processes
where each vertex has offspring distribution $D_n^\star-1$ and lifetime
distribution $Y$. By \cite{jagers-book}, the size of both flow
processes grow like $|\cF_i(t)| \sim\widetilde{W}_i^{(n)} \exp
(\lambda_n t)$ as $t\to\infty$, where $\lambda_n$ is the Malthusian
rate of growth of the branching process, and $\widetilde{W}_i^{(n)}
> 0$ ({due} to the fact that by assumption our branching processes
survive with probability $1$) are associated martingale limits.
Furthermore, an analysis of the two exploration processes suggests that
for $t>0$, the rate at which one flow cluster picks a vertex from the
other flow cluster (thus creating a collision in a small time interval
$[t,t+\mathrm{d}t)$) is approximately
%
\begin{equation}
\gamma_n(t) \approx\frac{\kappa|\cF_0(t)||\cF_1(t)|}{n} \approx \frac
{\kappa\widetilde{W}_0^{(n)} \widetilde{W}_1^{(n)} \exp
(2\lambda_n t)}{n},\qquad t
\geq0,
\end{equation}
where the constant $\kappa$ arises due to a subtle interaction of the
{stable-age} distribution of the associated continuous time
branching process with the exploration processes. This suggests that
times of creation of collision edge scales like $(2\lambda)^{-1} \log
{n}$, and further the time of birth of the first collision edge,
re-centered by $(2\lambda)^{-1}\log n$, converges to the first point
$\Xi_\infty$ of a Cox process with rate
\[
\gamma_\infty(x):= \kappa\widetilde{W}_0
\widetilde{W}_1 \exp( 2\lambda x),\qquad  x\in\bR.
\]
It is easy to check that
%
\begin{equation}
\label{eqXiinfty} \Xi_\infty\stackrel{d} {=} \frac{1}{2\lambda} (-\Lambda-
\log {\widetilde{W}_0} - \log{\widetilde{W}_1} + c ),
\end{equation}
where $c$ is a constant depending on $\lambda$ and $\kappa$, and
$\Lambda$ has Gumbel distribution independent of $\widetilde{W}_0,
\widetilde{W}_1$.

In \cite{BvdHHUnivFPP}, both $v_0$ and $v_1$ are chosen uniformly and
therefore have a degree \emph{distribution} different from the
{offspring} distribution $D_n^\star-1$ associated to the rest of the
branching process. Consequently, $\widetilde{W}_0^{(n)}$ and
$\widetilde{W}_1^{(n)}$ are not distributed as the martingale limit
$W^{(n)}$ but as a certain sum $\widehat{W}_s^{(n)}$ of such
variables (with $\widehat{W}_s^{(n)}\to\widehat{W}_s$ as $n\to
\infty
$). By contrast, in the setting of Proposition~\ref
{propWeightFinVarJointBiased} for $k=1$, the vertex $v_1$ has
distribution close to $D^\star-1$ by assumption, so that this
replacement is not necessary and $\widetilde{W}_1^{(n)}\stackrel{d}{=}
W^{(n)}\stackrel{d}{\longrightarrow}W$. Since the length of the
optimal path scales like
$2\Xi_n^{(1)}$, rearranging \eqref{eqXiinfty} gives Proposition~\ref{propWeightFinVarJointBiased} with $k=1$.

The actual rigorous proof in \cite{BvdHHUnivFPP} is a lot more subtle
albeit following the above underlying idea. The optimal path is formed
not quite at time $2\Xi_n$, and one has to keep track of ``residual
life-times'' of alive vertices, whose asymptotics follow from the
{stable-age-distribution} theory of Jagers and Nerman \cite
{jagers1984growth}, and so on, leading to the analysis of a much more
complicated Cox process. In the end, distributional identities for the
Poisson process yield the result above.

\subsubsection{Extension to multiple target vertices: The case \texorpdfstring{$k\geq2$}{$k>=2$}}

Let us now describe how one extends the above result for $k=1$ to more
general $k$. For ease of notation, assume $k=2$; the general case
follows in an identical fashion. Consider flow emanating from three
vertices $v_s$ and $v_1, v_2$ simultaneously at $t=0$. Arguing as
above, one finds that there exist paths $\cP_1$ and $\cP_2$ (not
necessarily optimal) between $v_s$ and $v_1, v_2$ such that the
respective lengths of the paths $\widetilde{\cC}_n ({v_{\widetildes s}},
v_1) $ and $\widetilde{\cC}_n ({v_{\widetildes s}}, v_2)$ satisfy
%
\begin{equation}\qquad
\label{eqntild-cn} \bigl(\lambda_n \widetilde{\cC}_n
({v_{\widetildes s}}, v_1) - \log {n} \bigr)_{i=1}^2
\stackrel{d} {\longrightarrow}(-\Lambda_i - \log {
\widehat{W}_s} - \log {W_i}+c)_{i=1}^2:=
\vW(2).
\end{equation}
Obviously the length of the optimal paths satisfy $C_n({v_{\widetildes
s}}, v_i) \leq\widetilde{\cC}_n ({v_{\widetildes s}},
v_1)$, and thus the limit $\vW(2)$ above serves as a limiting upper
bound (in the distributional sense) to the vector of lengths of optimal
costs properly re-centered,
\[
\vC_n(2):= \bigl(\lambda_n {\cC}_n
({v_{\widetildes s}}, v_1) - \log {n} \bigr)_{i=1}^2.
\]
However, the result holds for $k=1$ by the argument in the previous
section; thus the marginals of $\vC_n(2)$ must converge to the
marginals of $\vW(2)$. This implies that $\vC_n(2)$ converges to $\vW
(2)$. See \cite{Sale13} for more details.

\section{Convergence of the degree distribution}\label{spart-a-proofs}
In this section we prove part (a) of Theorems \ref
{tCMtheoremfinitevar}, \ref{tCMdegreesinfvar} and \ref
{tKnDegrees}, since the {proofs} share similarities. Parts (b) and~(c) of
these theorems are deferred to Sections~\ref{spart-b-proofs} and \ref{spart-c-proofs}. For the rest of the paper
we write
%
\begin{equation}
\label{eqnphiu-def} \phi_{W}(u):= \E\bigl(\exp(-uW)\bigr),\qquad u\geq0,
\end{equation}
for the Laplace transform of the random variable $W$ which arise as
martingale limits of branching processes and satisfy the recursive
distributional equations \eqref{eqWfors} or \eqref{eqWrecursionCMFinVar}.

All three proofs are based on an analysis of optimal path lengths,
using the following characterization of the out-degree of $v_t$. {Note
that we convert again to using $v_s$ for the source vertex, and we
continue to use $v_t$ for the target vertex.}

\begin{quote}
\textit{The out-degree of $v_t$ in $\TT_n$ is the number of
immediate neighbors of $v_t$ for which the shortest path from $v_s$
passes through vertex $v_t$.}
\end{quote}

To formalize this, write $\NN$ for the collection of neighbors of $v_t$
in $\GG_n$, and let $C_n'(v_s,v)$, $v\in\NN$, denote the shortest path
between vertices $v_s$ and $v$ in the modified graph $\GG'_n$ where the
vertex $v_t$, and all edges incident to $v_t$, are excised. Write
$Y_v$, $v\in\NN$, for the weight of the edge between $v$ and $v_t$; by
construction, the $Y_v$ are independent copies of $Y$, independent of
everything else. Then
%
\begin{equation}
\label{eqDistanceMinNeighbours} C_n(v_s,v_t)=\min
_{v\in\NN} \bigl( C_n'(v_s,v)+Y_v
\bigr),
\end{equation}
and the unique path in $\TT_n$ from $v_s$ to $v_t$ passes through the
unique vertex $U\in\NN$ for which $C_n(v_s,v_t)=C_n'(v_s,U)+Y_{
U}$. Moreover, the edge between $v_t$ and a vertex $v\in\NN\setminus
\{U\}$ belongs to $\TT_n$ if and only if the path from $v_s$ to $v$ via
$v_t$ is shorter than the optimal path excluding $v_t$. That is,
%
\begin{eqnarray}
\label{eqEdgeCondition} \mbox{edge }\{v_t, v\}&\in&\TT_n \quad\iff\quad
v=U \quad\mbox{or}
\nonumber
\\[-8pt]
\\[-8pt]
\nonumber
C_n'(v_s,U)+Y_{U}+Y_v&<&C_n'(v_s,v).
\end{eqnarray}
Because the alternatives in the right-hand side of \eqref
{eqEdgeCondition} are mutually exclusive, we can therefore express the
degree of $v_t$ as
%
\begin{equation}
\label{eqDegreevt} \deg_{\TT_n}(V_n) = 1+\sum
_{v\in\NN} \mathbh{1}_{\{C_n'(v_s,U)+Y_{U}+Y_v<C_n'(v_s,v)\}}.
\end{equation}

First we start with the configuration model. The proofs of part (a) of
Theorems~\ref{tCMtheoremfinitevar} and \ref{tCMdegreesinfvar} rely
on the asymptotics for optimal path lengths stated in Propositions~\ref
{propWeightFinVarJointBiased} and \ref{propWeightInfVarJointBiased}.

\begin{pf*}{Proof of Theorem~\ref{tCMtheoremfinitevar}(a)}
Since the original degree $d_{v_t}$ converges in distribution to $D$ as
$n\to\infty$, it suffices to condition on $\{d_{v_t}=k\}$ and then
show that $\deg_{\TT_n}(V_n)$ converges in distribution to $\widehat{D}_k$,
for each finite value $k\in\N$. Having made this conditioning, the event
%
\begin{equation}
A_{n,k}=\bigl\{d_{v_t}=k, v_t\neq
v_s, \vert\NN\vert=k, \NN\cap\{ v_s,v_t\}=
\varnothing\bigr\}
\end{equation}
(i.e., the event that the vertex paired to each of the $k$ stubs from
$v_t$, the vertex $v_s$, and the vertex $v_t$ itself are all distinct)
occurs with high probability.

It is easy to see that, conditional on the occurrence of $A_{n,k}$ and
the values $v_t$ and $\NN$, the graph $\GG'_n$ is equivalent to a
configuration model on the $n-1$ vertices $[n]\setminus\{v_t\}$,
where the degree $d'_v$ of vertex $v$ is given by
%
\begin{equation}
\label{eqDegreesWithoutvt} d'_v = \cases{ d_v-1,&\quad $v
\in\NN$,\vspace*{2pt}
\cr
d_v,&\quad $v\notin\NN$.}
\end{equation}
Conditional on $\{d_{v_t}=k\}$, let $v_1,\ldots,v_k$ denote the
vertices paired to stubs from~$v_t$. As discussed earlier, the vertices
$(v_1,\ldots,v_k)$ are chosen with probabilities asymptotically
proportional to $d_{v_1}\cdots d_{v_k}$. From \eqref
{eqDegreesWithoutvt} it follows that, conditional on $A_{n,k}$, the
modified degrees $(d'_{v_1},\ldots,d'_{v_k})$ converge jointly in
distribution to $k$ independent variables with the distribution
$D^\star
-1$; see \eqref{eqBiasedD}.
By Proposition~\ref{propWeightFinVarJointBiased}, conditional on
$A_{n,k}$, the recentered shortest paths $\lambda_{n-1}
C_n'(v_s,v_i)-\log(n-1)$, $i=1,\ldots,k$,
converge jointly in distribution to $-\log\widehat{W}_s-\log
W_i-\Lambda
_i+c$, $i=1,\ldots,k$,
while the edge weights $Y_{v_i}$ are independent copies of $Y$.
Recall the notation $M_k$ from \eqref{eqhatDk-CMFinVarDefn}. Then
\eqref{eqDistanceMinNeighbours} implies that
%
\begin{eqnarray}
\label{eqtocombine} %
&&\lambda_{n-1}
C_n(v_s,v_t)-\log(n-1)\nonumber\\
&&\qquad\stackrel{d} {
\longrightarrow }\min_{i=1,\ldots,k} (-\log\widehat{W}_s -
\log W_i -\Lambda_i+c+\lambda Y_i )
\\
&&\qquad=-M_k-\log\widehat{W}_s+c,\nonumber
\end{eqnarray}
also jointly with the previous convergences.

On the other hand, if we rescale and recenter the shortest paths in
\eqref{eqDegreevt}, then we get
%
\begin{eqnarray}
\label{eqDegreeScaled}&& \deg_{\TT_n}(V_n)
\nonumber
\\[-8pt]
\\[-8pt]
\nonumber
&&\qquad= 1+\sum
_{i=1}^k \mathbh{1}_{\{ ( \lambda_{n-1} C_n(v_s,v_t)-\log(n-1)  )
+\lambda_{n-1} Y_{v_i} <  ( \lambda_{n-1} C_n'(v_s,v_i) - \log
(n-1)  )\}}.
\end{eqnarray}
The mapping $(\lambda_{n-1} C_n'(v_s,v_i)-\log(n-1),Y_{v_i})_{i=1}^k
\mapsto\deg_{\TT_n}(V_n)$ defined by \eqref{eqDegreeScaled} is not
continuous. However, the limiting variables $(-\Lambda_i-\log\widehat
{W}_s-\log W_i+c,Y_i)_{i=1}^k$ are continuous and so is the difference
between the left and right-hand side of the inequality inside the
indicators; hence simple discontinuities of the mapping play no role.
By combining \eqref{eqtocombine} with \eqref{eqDegreeScaled}, as well
as the fact that $\lambda_n\rightarrow\lambda$, we conclude that,
conditional on $A_{n,k}$,
%
\begin{equation}
\deg_{\TT_n}(V_n) \stackrel{d} {\longrightarrow} 1+\sum
_{i=1}^k \mathbh{1}_{\{-M_k-\log\widehat{W}_s+c+\lambda Y_i <
-\Lambda_i-\log\widehat{W}_s-\log W_i+c\}},
\end{equation}
which simplifies to \eqref{eqhatDk-CMFinVarDefn}. Since
$d_{v_t}\stackrel{d}{\longrightarrow}
D$, this completes the proof of part (a).
\end{pf*}

Now we move to show the corresponding characterization of the degrees
in the shortest path tree in the {infinite-variance} case. The
proof is very similar, using Proposition~\ref
{propWeightInfVarJointBiased} in place of Proposition~\ref
{propWeightFinVarJointBiased}.

\begin{pf*}{Proof of Theorem~\ref{tCMdegreesinfvar}(a)}
For the infinite-variance case, no rescaling or recentering is needed
in \eqref{eqDegreevt}. Define $A_{n,k}$ and the modified shortest path
lengths $C'_n(v_s,v_i)$ as in the proof of Theorem~\ref
{tCMtheoremfinitevar}. Conditional on $A_{n,k}$, Proposition~\ref
{propWeightInfVarJointBiased} gives $(C'_n(v_s,v_i))_{i=1}^k\stackrel
{d}{\longrightarrow}
\widehat{V}_s+V_i$ and
%
\begin{equation}
C_n(v_s,v_t)\stackrel{d} {
\longrightarrow}\min_{i=1,\ldots,k} ( \widehat{V}_s+V_i+E_i
) = \widehat{V}_s+\xi_k,
\end{equation}
so that combining this with \eqref{eqDegreevt} gives that,
conditional on $A_{n,k}$,
%
\begin{equation}
\label{eqkDegreeLimitInfVar} \deg_{\TT_n}(V_n) \stackrel{d} {
\longrightarrow} 1+\sum_{i\neq U, 1\leq i\leq k} \mathbh{1}_{\{\widehat{V}_s+\xi
_k+E_i < \widehat{V}_s+V_i \}},
\end{equation}
which reduces to \eqref{eqhatDk-CMInfVarDefn} and completes the proof.
\end{pf*}

Now we aim to prove the similar characterization of the degrees for the
complete graph, that is, Theorem~\ref{tKnDegrees}(a).
The difficulty in this case is that the degree of $v_t$ is not tight,
and an additional argument is needed to show that only neighbors joined
to $v_t$ by short edges are likely to contribute to $\deg_{\TT_n}(V_n)$.

For the purposes of the following lemma, it is convenient to think of
$\TT_n$ as directed away from the source vertex $v_s$, so that the
children of $v_t$ are precisely those vertices $v$ for which $v_t$ is
the last vertex before $v$ on the shortest path from $v_s$ to $v$. In
this formulation, the out-degree of $v_t$ is equal to the number of
children of $v_t$ in $\TT_n$.

%
\begin{Lemma}\label{lemKnTightness}
Consider the complete graph with the edge cost distribution $E^s$, as
in Theorem~\ref{tKnDegrees}. Then, given $\varepsilon>0$, there exists
$R<\infty$ such that, with probability at least $1-\varepsilon$, every
edge between $v_t$ and a child of $v_t$ in the shortest-path tree $\TT
_n$ has edge weight at most $Rn^{-s}$.
\end{Lemma}

\begin{pf}
Let $\varepsilon>0$ be given. By Proposition~\ref{propWeightKnJoint}
applied for $k=1$, we may choose $R'<\infty$ such that $\log n -
R'\leq
\lambda_s n^s C_n(v_s,v_t)$ with probability at least $1-\frac
{1}{2}\varepsilon$. Assume {that} this event occurs, and suppose in
addition that $v_t$ has at least one child $V$ in $\TT_n$ joined to
$v_t$ by an edge of weights at least $Rn^{-s}$. Then
%
\begin{equation}
\lambda_s n^s C_n(v_s,V) \geq
\lambda_s n^s \bigl(C_n(v_s,v_t)+Rn^{-s}
\bigr) \geq\log n - R' + \lambda_s R,
\end{equation}
and furthermore $v_t$ is the last vertex before $V$ on the optimal
path from $v_s$ to $V$. Write $N$ for the number of vertices $v\in[n]$
with these two properties. Since $v_t$ is chosen uniformly,
independently of everything else,
%
\begin{equation}\qquad
\prob(N>0) \leq \Ev(N) \leq \sum_{v\in[n]}
\frac{1}{n} \prob\bigl(\lambda_s n^s
C_n(v_s,v) - \log n \geq \lambda_s R -
R'\bigr),
\end{equation}
and the right-hand side is the probability that a uniformly chosen
vertex $v$ has $\lambda_s n^s C_n(v_s,v) - \log n \geq\lambda_s R -
R'$. By Proposition~\ref{propWeightKnJoint} for $k=1$, this
probability can be made smaller than $\frac{1}{2}\varepsilon$ by
taking $R$ large enough.
\end{pf}

\begin{pf*}{Proof of Theorem~\ref{tKnDegrees}(a)}
For the collection of edges incident to $v_t$, write the edge weights
in increasing order as $E_1^s<\cdots<E_{n-1}^s$, and let $v_1,\ldots,v_{n-1}$ denote the corresponding ordering of the vertices
$[n]\setminus\{v_t\}$. It is easy to see that the rescaled order
statistics $(n-1)^s E_1^s, \ldots,(n-1)^s E_{n-1}^s$ converge to the
Poisson point process ${(X_i)_{i\geq1}}$ from \eqref{eqPPPEs}, in
the sense that for any $k\in\N$, jointly in $k$ and as $n\to\infty$,
%
\begin{equation}
\label{eqktuples-conv} \bigl((n-1)^s E_1^s,
\ldots,(n-1)^s E_k^s \bigr) \stackrel{d}{\longrightarrow}(X_1,\ldots,X_k).
\end{equation}
This follows from the usual convergence of the rescaled order
statistics $(n-1)E_1<\cdots<(n-1)E_{n-1}$ toward a Poisson point process
of unit intensity, together with the fact that the map $x\mapsto x^s$
is increasing and continuous.

If $v_t$ had only a fixed number $k$ of neighbors, we could complete
the proof in the same way as for Theorems~\ref{tCMtheoremfinitevar}
and \ref{tCMdegreesinfvar}. We must therefore control the
possibilities that (a)~some vertex not belonging to $\{v_1,\ldots,v_k\}$ (for some $k$) contributes to the out-degree of $v_t$; and
(b)~the last vertex before $v_t$ on the shortest path from $v_s$ to
$v_t$ does not belong to $\{v_1,\ldots,v_k\}$ for some $k$.

Let $B_{n,k}$ denote the event that every child of $v_t$ in $\TT_n$ is
one of the vertices $\{v_1,\ldots,v_k\}$. We claim that
%
\begin{equation}
\label{eqBnkBound} \lim_{k\to\infty} \liminf_{n\to\infty}
\prob(B_{n,k})=1.
\end{equation}
Indeed, by a union bound we have that if $B_{n,k}^c$ occurs, then
either the $k$th edge weight is too small, or if it is not, then $v_t$
has a neighbor in $\TT_n$ with too large edge-weight
%
\begin{equation}
\label{equnionbound} \Pv\bigl(B_{n,k}^c\bigr) \le\Pv\bigl(
n^s E_k^s\leq R \bigr) + \Pv\bigl(
v_t \mbox{ has a child } v \mbox{ with } Y_v \ge R
n^{-s} \bigr).
\end{equation}
But from \eqref{eqktuples-conv} we know that $n^s E_k^s \stackrel
{d}{\longrightarrow}X_k$ as
$n\to\infty$ (the distinction between $n$ and $n-1$ being irrelevant in
this limit). Since $X_k\stackrel{{\mathbb P}}{\longrightarrow}\infty
$ as $k\to\infty$, we can choose
$R=R(k)$ in such a way that
%
\begin{equation}
\lim_{k\to\infty}\limsup_{n\to\infty}\prob
\bigl(n^s E_k^s \leq R(k)\bigr) =0,
\end{equation}
and then Lemma~\ref{lemKnTightness} shows that the second term in
\eqref{equnionbound} is also negligible; hence we get \eqref{eqBnkBound}.

On $B_{n,k}$, only the vertices $v_1,\ldots,v_k$ contribute to the
out-degree of $v_t$, and \eqref{eqDegreevt} becomes
%
\begin{equation}
\label{eqOnAnk} \deg_{\TT_n}(V_n)= 1 + \sum
_{i=1}^k \mathbh{1}_{\{
C_n(v_s,v_t)+E_i^s<C'_n(v_s,v_i)\}} \qquad\mbox{on
$B_{n,k}$}.
\end{equation}

Since the original graph is the complete graph, the modified graph $\GG
'_n$ with $v_t$ excluded is a complete graph on the $n-1$ vertices
$[n]\setminus\{v_t\}$. Since the labeling of $v_1,\ldots,v_k$ depend
only on the excluded edge weights, Proposition~\ref{propWeightKnJoint}
applies, and we conclude that
%
\begin{eqnarray}
\label{eqKnRescaledktuple} &&\bigl(\lambda_s (n-1)^s
C'_n(v_s,v_i)-\log(n-1),
(n-1)^s E_i^s\bigr)_{i=1}^k
\nonumber
\\[-8pt]
\\[-8pt]
\nonumber
&&\qquad
\stackrel{d} {\longrightarrow}(-\Lambda_i-\log W_s-\log
W_i+c, X_i)_{i=1}^k.
\end{eqnarray}
We wish to conclude also that
%
\begin{equation}
\label{eqCnConvergesKn} \lambda_s (n-1)^s C_n(v_s,v_t)-
\log(n-1) \stackrel {d} {\longrightarrow}-M-\log W_s+c,
\end{equation}
jointly with the convergence in \eqref{eqKnRescaledktuple}. However,
\eqref{eqCnConvergesKn} does not follow from \eqref
{eqDistanceMinNeighbours} and \eqref{eqKnRescaledktuple}; rather, we
obtain only that
%
\begin{eqnarray}
\label{eqCnTruncated}&& \lambda_s (n-1)^s \min
_{i=1,\ldots,k} \bigl( C'_n(v_s,v_i)+E_i^s
\bigr)
\nonumber
\\[-8pt]
\\[-8pt]
\nonumber
&&\qquad\stackrel{d} {\longrightarrow}-\max_{i=1,\ldots,k} ( \Lambda
_i+\log W_i-\lambda_s X_i )-\log
W_s+c,
\end{eqnarray}
that is, the maximum is taken only on the first $k$ elements.
We will therefore give a separate argument to show \eqref{eqCnConvergesKn}.

Set $M'_k=\max_{{i\in[k]}} ( \Lambda_i+\log W_i-\lambda_s X_i
 )$, so that $M=\sup_k M'_k$. Further, let $( Z,(-\Lambda
_i-\log
W_s-\log W_i+c, X_i)_{i\geq1})$ denote any subsequential limit of the
rescaled shortest paths
%
\begin{eqnarray}
&&\bigl(\lambda_s (n-1)^s C_n(v_s,v_t)-
\log(n-1),
\nonumber
\\[-8pt]
\\[-8pt]
\nonumber
&&\qquad \bigl(\lambda_s (n-1)^s C'_n(v_s,v_i)-
\log(n-1), (n-1)^s E_i^s \bigr)_{{i\in[n-1]}}
\bigr).
\end{eqnarray}
By \eqref{eqCnTruncated}, $Z\leq-M'_k-\log W_s+c$ for each $k$, and
therefore $Z\leq-M-\log W_s+c$. It therefore suffices to show that the
\emph{marginal distribution} of $Z$ is the same as that of $-M-\log
W_s+c$. The event that $M<m$ is the event that the number of points
$(X_i,\Lambda_i+\log W_i)$ lying in the region $\{(x,y)\colon
y-\lambda_s x \geq m\}$ should be $0$. Since $(\Lambda_i)_{i\geq
1},(W_i)_{i\geq1}$ are i.i.d., the collection $(X_i,\Lambda_i+\log
W_i)_{i\geq1}$ forms a Poisson point process on $(0,\infty)^2$ with
intensity measure $\mathrm{ d}\mu_s\times\prob(\Lambda+\log W\in\cdot)$,
and we compute
%
\begin{eqnarray}
\label{eqM-logW-DF} %
\prob(M<m) &=& \exp \biggl( -\int
_0^\infty\prob(\Lambda+\log W\geq
\lambda_s x+m)\,\mathrm{ d}\mu_s(x) \biggr)
\nonumber\\
& =& \exp \biggl( - \int_0^\infty\prob(\log E
\leq-\lambda_s x-m+\log W)\,\mathrm{ d}\mu_s(x) \biggr)
\nonumber
\\[-8pt]
\\[-8pt]
\nonumber
& = &\exp \biggl( -\int_0^\infty\E \bigl( 1-\exp
\bigl( - W{\mathrm e}^{-\lambda_s x-m} \bigr) \bigr)\,\mathrm{ d}\mu_s(x)
\biggr)
\\
& = &\exp \biggl( -\int_0^\infty \bigl( 1-
\phi_{W} \bigl({\mathrm e}^{-\lambda_s
x-m} \bigr) \bigr) \,\mathrm{ d}
\mu_s(x) \biggr),\nonumber
\end{eqnarray}
where $\phi_{W}(u)=\E({\mathrm e}^{-uW})$ is the Laplace
transform of $W$.
The recursive definition~\eqref{eqWfors} of $W$ implies the identity
%
\begin{equation}
\label{eqphiWrecursion} \phi_{W}(u) = \exp \biggl( -\int_0^\infty
\bigl( 1-\phi_{
W}\bigl(u{\mathrm e} ^{-\lambda_s x}\bigr) \bigr)\,\mathrm{
d}\mu_s(x) \biggr),
\end{equation}
so that \eqref{eqM-logW-DF} reduces to
%
\begin{equation}
\label{eqMdistribution} \prob(M<m)=\phi_{W}\bigl({\mathrm e}^{-m}
\bigr).
\end{equation}
In particular, we have $\prob(-M-\log W_s+c > z)=\E ( \phi
_{
W}(W_s {\mathrm e}^{z-c})  )$.

On the other hand, since $Z$ is the limit in distribution of $\lambda_s
(n-1)^s C_n(v_s,v_t)-\log n$, Proposition~\ref{propWeightKnJoint}
implies that $Z \stackrel{d}{=}-\Lambda-\log W_s-\log W_t+c$ (the distinction
between $n$ and $n-1$ again being irrelevant), and we compute
%
\begin{eqnarray}
\label{eqZdistribution} \prob(Z > z) &=& \E \bigl( \mathbb{P}(\log E-\log W_s-
\log W > z-c \vert W_s,W) \bigr)
\nonumber
\\[-8pt]
\\[-8pt]
\nonumber
& = &\E \bigl( \exp \bigl( -W_s W {\mathrm e}^{z-c} \bigr)
\bigr) =\E \bigl( \phi_{W}\bigl(W_s {\mathrm
e}^{z-c}\bigr) \bigr).
\end{eqnarray}
This proves \eqref{eqCnConvergesKn}.

We can now complete the proof of Theorem~\ref{tKnDegrees}(a). Rescale
and recenter the edge weights, and apply \eqref
{eqKnRescaledktuple} and \eqref{eqCnConvergesKn} to the right-hand side
of \eqref{eqOnAnk} to conclude that, on $B_{n,k}$, $\deg_{\TT
_n}(V_n)$ is
equal to a random variable that converges in distribution to
%
\begin{equation}
\widetilde{D}_k = 1+\sum_{i=1}^k
\mathbh{1}_{\{\Lambda_i+\log
W_i+\lambda_s X_i < M\}}.
\end{equation}
Since $\widehat{D}$ is finite a.s., $\prob(\widetilde{D}_k\neq
\widehat
{D})\to0$ as $k\to\infty$. Together with \eqref{eqBnkBound}, this
completes the proof.
\end{pf*}

In the course of proving \eqref{eqCnConvergesKn} [compare \eqref
{eqMdistribution} with the calculation in \eqref{eqZdistribution}],
we have proved an equality in law between $M$ and $\Lambda+\log W$,
which we record for future reference:

\begin{Lemma}\label{lemMLambdalogW}
The random variables $M$ and $W$ from Section~\ref{ssKnResults} are
related by
%
\begin{equation}
M \stackrel{d} {=}\Lambda+ \log W.
\end{equation}
\end{Lemma}

Observe that the result of Lemma~\ref{lemMLambdalogW} does not apply
in the CM setting from Section~\ref{ssCMfiniteVarResults} because of
size-biasing and depletion-of-points effects.

\section{Convergence of the empirical degree distribution}\label
{spart-b-proofs}
In this section we sketch the proofs of part (b) of Theorems \ref
{tKnDegrees}, \ref{tCMtheoremfinitevar} and \ref
{tCMdegreesinfvar}. Since $v_t$ is a uniformly chosen vertex, $\E
(\widehat{p}_k^{(n)})=\prob(\deg_{\TT_n}(V_n)=k)\to\prob(\widehat{D}=k)$
by part (a). By an application of Chebychev's inequality, it suffices
to prove that
%
\begin{equation}
\label{eqTwoDegrees} \prob\bigl(\deg_{\TT_n}(v_t)=k,
\deg_{\TT_n}(w_t)=k\bigr) \to\prob(\widehat
{D}=k)^2,
\end{equation}
where $w_t$ is another uniformly chosen vertex independent of $v_t$.

\begin{pf*}{Proof of Theorem~\ref{tCMtheoremfinitevar}(b)}
As in the proof of part (a), it suffices to condition on the original
degrees. Fix $i,j\in\N$. Conditional on $\{d_{v_t}=i,d_{w_t}=j\}$,
the event
\[
A_{n, i,j}= \bigl\{ d_{v_t}=i, d_{w_t}=j,
v_t, w_t, v_s, \CN(v_t),
\CN(w_t) \mbox{ all distinct}\bigr\}
\]
occurs with high probability. Moreover, Proposition~\ref
{propWeightFinVarJointBiased} holds for the $i+j$ neighbors of $v_t$
and $w_t$, saying that the re-centered edge weights tend to
exchangeable random variables.
As in \eqref{eqtocombine} and \eqref{eqDegreeScaled}, we get that,
conditionally on $A_{n,i,j}$,
%
\begin{eqnarray}
\label{eqjointconv} %
 \deg_{\TT_n}(v_t)
&\stackrel{d}{\longrightarrow}&1+\sum_{l=1}^i
\mathbh{1}_{\{-M_i^{(v_t)} -\log\widehat W_s +c +\lambda_n
Y_l < - \La_l - \log W_l -\log\widehat W_s +c\}},
\nonumber
\\[-8pt]
\\[-8pt]
\nonumber
\deg_{\TT_n}(w_t) &\stackrel{d}{\longrightarrow}&1+
\sum_{l=i+1}^{i+j} \mathbh{1}_{\{-M_j^{(w_t)} -\log\widehat W_s +c
+\lambda_n Y_l < - \La_l - \log W_l -\log\widehat W_s +c\}},
\end{eqnarray}
where $M_i^{(v_t)}= \max_{l=1, \ldots, i} (\La_l + \log W_l - \lambda Y_l)$
and $M_j^{(w_t)}= \max_{l=i+1, \ldots, i+j} (\La_l + \log W_l -
\lambda
Y_l)$. The terms $\log\widehat{W}_s$ cancel in \eqref{eqjointconv},
and it follows that $\deg_{\TT_n}(V_n)$ and $\deg_{\TT_n}(w_t)$
converge to
independent limits conditional on $\{d_{v_t}=i,d_{w_t}=j\}$. By
Condition \ref{condCMFinVar}, and since $v_t$ and $w_t$ are both
independent uniform draws from~$[n]$, the random variables $d_{v_t}$
and $d_{w_t}$ converge jointly to independent copies of $D$. Thus it
follows that $\deg_{\TT_n}(v_t)$ and $\deg_{\TT_n}(w_t)$ converge
(unconditionally) to independent copies of~$\widehat{D}$. In
particular, \eqref{eqTwoDegrees} holds.
\end{pf*}

The proof of Theorem~\ref{tCMdegreesinfvar}(b) is identical, using
Proposition~\ref{propWeightInfVarJointBiased} instead of Proposition~\ref{propWeightFinVarJointBiased} as in the proof of part (a):

\begin{pf*}{Proof of Theorem~\ref{tKnDegrees}(b)}
The idea here is again similar to the proof of Theorem~\ref
{tKnDegrees}(a). First, arrange the outgoing edge weights from $v_t$
and $w_t$ separately in increasing order and multiply by $(n-2)^s$.
Since the weight of the edge between $v_t$ and $w_t$ diverges under
this rescaling, we see that these rescaled edge weights converge to two
\emph{independent} Poisson processes ${(X_i^{(v_t)})_{i\geq1}}$
and ${(X_i^{(w_t)})_{i\geq1}}$. Denote the corresponding two
orderings of vertices by ${(v_i)_{i\geq1}}$ and ${(w_i)_{i\geq1}}$.
For any fixed $k\in\N$, the vertices $v_s,v_t,w_t,v_1,\ldots,v_k,w_1,\ldots,w_k$ are all distinct with high probability, and
conditional on this event we can apply Proposition~\ref
{propWeightKnJoint} to the $2k$ vertices $v_1,\ldots,v_k,w_1,\ldots,w_k$. A modification of the argument from the proof of part (a), as in
the discussion following \eqref{eqCnConvergesKn}, shows that $\lambda
_s (n-2)^s C_n(v_s,v_t)-\log(n-2)$ and $\lambda_s (n-2)^s
C_n(v_s,w_t)-\log(n-2)$ converge jointly to $-M^{(v_t)}-\log W_s+c$ and
$-M^{(w_t)}-\log W_s+c$, where $M^{(v_t)}$ and $M^{(w_t)}$ are
independent; we leave the details to the reader. With $B_{n,k}^{(w_t)}$
denoting the analogue of $B_{n,k}$ [where $B_{n,k}$ is defined above
\eqref{eqBnkBound}] with $v_t$ replaced by $w_t$, we conclude from
\eqref{eqOnAnk} that, on $B_{n,k}\cap B_{n,k}^{(w_t)}$, $\deg_{\TT_n}(V_n)$
and $\deg_{\TT_n}(w_t)$ are equal to random variables that converge in
distribution to independent copies of $\widetilde{D}_k$. Since
$B_{n,k}$ and $B_{n,k}^{(w_t)}$ both satisfy~\eqref{eqBnkBound}, we
conclude that $\deg_{\TT_n}(V_n)$, $\deg_{\TT_n}(w_t)$ have independent
limits, and~\eqref{eqTwoDegrees} holds.
\end{pf*}

\section{Average degrees}\label{spart-c-proofs}
In this section we prove part (c) of Theorems \ref{tKnDegrees}, \ref
{tCMtheoremfinitevar} and \ref{tCMdegreesinfvar}.
Here we show that the average of the limiting degree in all the three
cases is $2$, as one would expect.

\begin{pf*}{Proof of Theorem~\ref{tKnDegrees}(c)}
Recall that $\mu_s$ stands for the intensity measure for the ordered
points $X_i$, as in Section~\ref{ssKnResults}, and recall the
characterization of the degree $\widehat D$ in part (a) of Theorem~\ref
{tKnDegrees}. Since ${(\Lambda_i+\log W_i)_{i\geq1}}$ are
i.i.d. random variables, the points ${(X_i, \Lambda_i+\log W_i)_{i\geq1}}$
form a Poisson point process (PPP) $\cP$ on $\R^+\times\R$ with
the product intensity measure $\mu_s(\dd x)\cdot\Pv (\Lambda
+\log
W\in\dd y  )$; see, for example, \cite{Resnick86}, Proposition~2.2.

The event that $\{M> m\}$ is the event that the number of points
$X_i,\Lambda_i+\log W_i$ lying in the region $\{(x',y')\colon
y'-\lambda_s x' > m\}$ is at least $1$. Hence $\{M> m\}$ is measurable
with respect to the $\sigma$-field generated by the restriction of the
Poisson point process ${(\Lambda_i+\log W_i)_{i\geq1}}$ to the
infinite upward-facing triangle\vspace*{1pt} $\{(x',y') \in\R^+ \times\R\dvtx y'-\lambda_s x' > m \}$.
On the other hand, a point $(x,y)$ contributes to $\widehat D$ if $M >
y+\lambda_s x$, and clearly the point $(x,y)$ does not lie in the
infinite upward-facing triangle
$\{(x',y') \in\R^+ \times\R\dvtx y'-\lambda_s x' > m=y+\lambda_s x
\}$.

Hence, by the independence of PPP points in disjoint sets, conditional
on finding a point $(X,\Lambda+\log W)$ with value $(x,y)$, the
conditional probability of the event $\{M > y+\lambda_s x\}$ is equal
to the unconditional probability, which is $1-\phi_{W}({\mathrm e}
^{-y-\lambda_s x})$ by Lemma~\ref{lemMLambdalogW}.
On the other hand, $\prob(\Lambda+\log W\leq y)=\phi_{
W}({\mathrm e}
^{-y})$ implies that the intensity measure for the points $(X,\Lambda
+\log W)$ is $\mathrm{ d}\mu_s(x)\times(-\phi_{W}'({\mathrm
e}^{-y})){\mathrm e}
^{-y}\,\mathrm{ d}y$. Hence
%
\begin{eqnarray}
\E(\widehat{D}-1) &=& \E \biggl( \sum_{(x,y)\in\{(X_i,\Lambda_i+\log W_i), i\in\N\}} \prob(M
> y+\lambda_s x) \biggr)
\nonumber
\\
& =& \int_0^\infty\int_{-\infty}^\infty
\bigl(1-\phi_{W}\bigl({\mathrm e}^{-y-\lambda_s
x}\bigr)\bigr) \bigl(-
\phi_{W}'\bigl({\mathrm e}^{-y}\bigr)\bigr){
\mathrm e}^{-y}\,\mathrm{ d}y \,\mathrm{ d}\mu_s(x)
\\
& = &\int_0^\infty\int_0^\infty
\bigl(1-\phi_{W}\bigl(u{\mathrm e}^{-\lambda
_s x}\bigr)\bigr) \bigl(-
\phi _{W}'(u)\bigr)\,\mathrm{ d}u \,\mathrm{ d}\mu_s(x)\nonumber
\end{eqnarray}
by the substitution $u={\mathrm e}^{-y}$. By relation \eqref{eqphiWrecursion},
we obtain
\[
\E(\widehat{D}-1) = \int_0^\infty\bigl(-\log
\bigl(\phi_{W}(u)\bigr)\bigr) \bigl(-\phi_{W}'(u)
\bigr)\,\mathrm{ d}u = \int_0^1 (-\log x)\, \mathrm{ d}x
=1. 
\]
\upqed\end{pf*}

Next we give a direct proof of the average degree in shortest path tree
for the configuration model with finite-variance degrees.

\begin{pf*}{Proof of Theorem~\ref{tCMtheoremfinitevar}(c)}
Let $f(z)=\E(z^D)$ denote the probability generating function of $D$.
Then the probability generating function of $D^\star-1$ is
$f'(z)/f'(1)$, and from \eqref{eqWrecursionCMFinVar} it follows that
%
\begin{equation}
\label{eqphiWrecursionCMFinVar} \phi_{W}(u)=\frac{f' ( \E(\phi_{W}(u{\mathrm
e}^{-\lambda Y}))
 )}{f'(1)}.
\end{equation}
In \eqref{eqhatD-CMFinVarDefn}, partition according to the value of
$D$ and use symmetry to see that
%
\begin{eqnarray}
&&\E(\widehat{D}-1)\nonumber\\
&&\qquad = \sum_{k=1}^\infty
\prob(D=k) k\prob(\Lambda_1+\log W_1+\lambda
Y_1 < M_k)
\nonumber
\\
&&\qquad = \sum_{k=1}^\infty\prob(D=k)\nonumber\\
&&\hspace*{46pt}{}\times k \bigl( 1-
\prob(\Lambda_i-\log W_i-\lambda Y_i \leq
\Lambda_1+\log W_1+\lambda Y_1, i=2,\ldots,k)
\bigr)
\nonumber
\\[-8pt]
\\[-8pt]
\nonumber
&&\qquad = \sum_{k=1}^\infty\prob(D=k)
\\
&&\hspace*{46pt}{}\times k\E \bigl(
1-\mathbb{P}(\Lambda+\log W-\lambda Y\leq\Lambda_1+\log
W_1+\lambda Y_1 \vert\Lambda _1,W_1,Y_1)^{k-1}
\bigr)
\nonumber\\
&&\qquad = \E \bigl( f'(1)-f' \bigl( \mathbb{E}\bigl(
\phi_{W}\bigl({\mathrm e}^{-\lambda
Y-\Lambda_1-\log W_1-\lambda Y_1}\bigr) \vert
\Lambda_1,W_1,Y_1\bigr) \bigr) \bigr)
\nonumber
\\
&&\qquad = f'(1)\E \bigl( 1-\phi_{W}\bigl({\mathrm
e}^{-\Lambda_1-\log
W_1-\lambda Y_1}\bigr) \bigr),\nonumber
\end{eqnarray}
by \eqref{eqphiWrecursionCMFinVar} with $u={\mathrm e}^{-\Lambda
_1-\log
W_1-\lambda Y_1}$. Integrating first over $Y_1$ and using $\prob
(\Lambda
+\log W\leq x)=\phi_{W}({\mathrm e}^{-x})$ and \eqref
{eqphiWrecursionCMFinVar} again,
%
\begin{eqnarray}
\E(\widehat{D}-1) &=& f'(1)\int_{-\infty}^\infty
\bigl( 1-\E \bigl( \phi_{
W}\bigl({\mathrm e} ^{-x-\lambda Y}\bigr)
\bigr) \bigr) \bigl(-\phi_{W}'\bigl({\mathrm
e}^{-x}\bigr){\mathrm e}^{-x}\bigr) \,\mathrm{ d}x
\nonumber
\\
& =& f'(1)\int_0^\infty \bigl( 1-\E
\bigl( \phi_{W}\bigl(u {\mathrm e}^{-\lambda Y}\bigr) \bigr) \bigr)
\bigl(-\phi_{W}'(u)\bigr) \,\mathrm{ d}u
\nonumber
\\[-8pt]
\\[-8pt]
\nonumber
& =& f'(1)\int_0^\infty \bigl( 1-
\bigl(f' \bigr)^{-1} \bigl( f'(1)\phi
_{
W}(u) \bigr) \bigr) \bigl(-\phi_{W}'(u)
\bigr) \,\mathrm{ d}u
\\
& = &\int_0^1 (1-z) f''(z)
\,\mathrm{ d}z = \bigl[ (1-z)f'(z) \bigr]_0^1 +
\int_0^1 f'(z) =1,\nonumber
\end{eqnarray}
where we used the substitution $f'(z)=f'(1)\phi_{W}(u)$. Here we
used $f'(0)=f(0)=0$, which follows from the assumption that $D\geq2$ a.s.
\end{pf*}

Next we give a direct proof for the average degree in the shortest path
tree for the configuration model with infinite-variance degrees.

\begin{pf*}{Proof of Theorem~\ref{tCMdegreesinfvar}(c)}
In the setting of Theorem~\ref{tCMdegreesinfvar} it is relevant to
consider the distribution function
$F_{V}(x)=\prob(V\leq x)$ instead of the Laplace transform. Then
from \eqref{eqexplosiontime-representation} we obtain
%
\begin{equation}
\label{eqFVrecursion} 1-F_{V}(x)=\frac{f'(\P(V+E>x))}{f'(1)}.
\end{equation}
Partition \eqref{eqhatD-CMInfVarDefn} according to the value of $D$,
and use the continuity of the distibutions to obtain
%
\begin{eqnarray}
\E(\widehat{D}-1) &=& \sum_{k=1}^\infty
\prob(D=k) k\prob(V_1-E_1 > \xi_k)
\nonumber
\\
& = &\sum_{k=1}^\infty\prob(D=k) k\E \bigl(
1-\mathbb{P}(V_i+E_i \leq V_1-E_1
\vert V_1,E_1)^{k-1} \bigr)
\nonumber
\\[-8pt]
\\[-8pt]
\nonumber
& =& \E \bigl( f'(1)-f' \bigl( \mathbb{P}(V+E \leq
V_1-E_1 \vert V_1,E_1) \bigr)
\bigr)
\\
& = &f'(1)\E \bigl( F_{V} ( V_1-E_1
) \bigr),\nonumber
\end{eqnarray}
where we applied \eqref{eqFVrecursion} with $x=V_1-E_1$. That is,
using \eqref{eqFVrecursion} again,
%
\begin{eqnarray}
\E(\widehat{D}-1) &= &f'(1)\prob(V\leq V_1-E_1)=f'(1)
\prob(V+E_1\leq V_1)
\nonumber
\\
& =& f'(1) \int_0^\infty\prob(V+E
\leq x) F_{V}'(x)\,\dd x
\nonumber
\\[-8pt]
\\[-8pt]
\nonumber
& =& f'(1) \int_0^\infty \bigl[
\bigl(f' \bigr)^{-1} \bigl( f'(1)
\bigl(1-F_{
V}(x)\bigr) \bigr) \bigr] F_{V}'(x)
\,\dd x
\\
& =&\int_0^1 (1-z)f''(z)
\,\dd z=1\nonumber
\end{eqnarray}
as before, where we used the substitution $f'(z)=f'(1)(1-F_{V}(x))$.
\end{pf*}

\begin{Remark}
An alternative proof of part (c) of Theorems \ref
{tCMtheoremfinitevar} and \ref{tCMdegreesinfvar} is the following:
Because $v_t$ is a uniformly chosen vertex, we have
\[
\Ev_n\bigl(\deg_{\TT_n}(V_n)\bigr) = \Ev
\biggl( \frac{1}{n} \sum_{v\in[n]} \deg
_{\TT
_n}(v) \biggr).
\]
The sum of the degrees is twice the number of edges, namely $2(n-1)$
since $\TT_n$ is a tree on $n$ vertices. Therefore $\E(\deg_{\TT_n}(V_n)
)\to
2$. On the other hand, we have $\deg_{\TT_n}(V_n)\stackrel{d}{\longrightarrow}\widehat{D}$ and
$\deg_{\TT_n}(V_n)\leq d_{v_t}\stackrel{d}{\longrightarrow
}D$. Under the hypotheses of Theorem~\ref
{tCMtheoremfinitevar} or Theorem~\ref{tCMdegreesinfvar}, $D$ has
finite expectation, and we can make a dominated convergence argument to
show that $\E(\deg_{\TT_n}(V_n))\to\E(\widehat{D})$.
Note that this reasoning is not available on the complete graph, where
the original degree $d_{v_t}$ diverges.
\end{Remark}

\section{Degree asymptotics}\label{ssasymptotics}
In this section we prove the theorems investigating the asymptotic
behavior of the degrees in the shortest path tree.
\subsection{Degree asymptotics: CM with finite-variance degrees}

Now we prove Theorems \ref{tCMdegreekfin} and \ref{tRateOfConv}.
{Theorem~\ref{tCMdegreekfin}} tells us that almost all the edges
of a large degree vertex are revealed by the shortest path tree.
{Theorem~\ref{tRateOfConv}} shows that the finite order correction
term, that is, the number of ``hidden'' edges, still can be quite large
{for some edge-weight distributions.} The main advantage is that in
both cases we can use the representation of the degrees in Theorem~\ref
{tCMtheoremfinitevar}(a).

\begin{pf*}{Proof of Theorem~\ref{tCMdegreekfin}}
We have $\prob(\Lambda+\log W > x)>0$ for each $x\in\Rbold$, by either
of the hypotheses on $\Lambda$ or $W$. It follows that $M_k\stackrel
{{\mathbb P}}{\longrightarrow}
\infty
$ as $k\to\infty$. Let $\varepsilon>0$ be given, and choose
$x<\infty$
such that $q=\prob(\Lambda+\log W+\lambda Y<x)$ satisfies $q\geq
1-\varepsilon$. Then
%
\begin{equation}
\label{eqDstarBinomialLower} \widehat{D}_k\geq\sum_{i=1}^k
\mathbh{1}_{\{\Lambda_i+\log
W_i+\lambda Y_i<x\}}\qquad \mbox{on }\{M_k>x\},
\end{equation}
and the right-hand side of inequality \eqref{eqDstarBinomialLower} is
$\operatorname{Binomial}(k,q)$. Since $\prob(M_k>x)\to1$, it follows that $\prob
(\widehat{D}_k\geq k(1-2\varepsilon))\to1$, and since $\varepsilon>0$
was arbitrary, this shows that $\widehat{D}_k=k(1-o_{\prob}(1))$.
\end{pf*}

\begin{pf*}{Proof of Theorem~\ref{tRateOfConv}}
For part (a), recall that $M_k$ is the maximum of~$k$
i.i.d. random variables $\Lambda_i+\log W_i-\lambda Y_i$, so, by
classical extreme value theory~\cite{EmbKluMik97}, $M_k=\log
k+O_{\prob}(1)$
will follow if $c{\mathrm e}^{-x} \leq\prob(\Lambda+\log W-\lambda
Y> x)\leq C
{\mathrm e}^{-x}$ for $x$ sufficiently large. For the upper bound,
write $\Lambda
=-\log E$ and use $\prob(E<x)\leq x$ for $x>0$ to obtain
%
\begin{eqnarray}
 \prob(\Lambda+\log W-\lambda Y>x)&=&\E \bigl(
\mathbb{P}\bigl(E< \smash{W{\mathrm e} ^{-\lambda Y}{\mathrm e}^{-x}}
\vert W,Y\bigr) \bigr)
\nonumber
\\[-8pt]
\\[-8pt]
\nonumber
&\leq&\E \bigl( \smash{W{\mathrm e}^{-\lambda Y}{\mathrm e}^{-x}}
\bigr)=O\bigl({\mathrm e}^{-x}\bigr).
\end{eqnarray}
The lower bound follows from $\prob(E<y)\geq cy$ for some $c>0$,
uniformly over $y<1$:
%
\begin{equation}
\prob(\Lambda+\log W-\lambda Y>x) \geq\E \bigl( \mathbh{1}_{\{W<K\}}
\smash{cW{\mathrm e}^{-\lambda
Y}{\mathrm e}^{-x}} \bigr) \geq
c' {\mathrm e}^{-x}
\end{equation}
for $K$ large enough and $x$ large enough that $K{\mathrm e}^{-x}\leq
1$. This
completes the proof of part (a).

For part (b), let $\varepsilon>0$ be given,
and choose $K<\infty$ large enough that $\prob(M_k<\log
k-K)<\varepsilon
$. Apply \eqref{eqDstarBinomialLower} with $x=\log k-K$ to conclude
that, apart from an event of small probability, $\widehat{D}_k$ is
stochastically larger than a $\operatorname{Binomial}(k,p_k)$ random variable with
$p_k=\prob(\Lambda+\log W+\lambda Y < \log k-K)$. To show tightness for
$k-\widehat{D}_k$, it is therefore sufficient to show that
$1-p_k=O(1/k)$. [To see the sufficiency, note that we need only show
that the $\operatorname{Binomial}(k,1-p_k)$ distributions are tight, and
$1-p_k=O(1/k)$ implies that these distributions have a uniformly
bounded mean. Alternatively, note that the $\operatorname{Binomial}(k,C/k)$
distribution converges to the Poisson$(C)$ distribution as $k\to\infty
$.] We compute
%
\begin{eqnarray}
1-p_k &=& \E \bigl( \mathbb{P}(-\Lambda
\leq\log W+\lambda Y-\log k+K \vert W,Y) \bigr)
\nonumber\\
&=& \E \biggl( \mathbb{P}\biggl(E\leq \smash{\frac{1}{k}W{\mathrm
e}^{\lambda Y}{\mathrm e} ^K} \Big\vert W,Y\biggr) \biggr)
\\
& \leq& O\bigl(k^{-1}\bigr)\E \bigl( \smash{W{\mathrm e}^{\lambda Y}}
\bigr)=O\bigl(k^{-1}\bigr), \nonumber
\end{eqnarray}
since $\Ev(W {\mathrm e}^{\lambda Y})<\infty$ by assumption.

For part (c), suppose $\lambda>1$. For the upper
bound, we estimate
%
\begin{eqnarray}
1-p_k& =& \E \bigl( \mathbb{P}(\lambda Y \geq\log k-K-\Lambda-\log W
\vert \Lambda,W) \bigr)
\nonumber
\\
& =& \prob(\log k-K-\Lambda-\log W<0)
\nonumber
\\
& &{}+\E \bigl( \mathbh{1}_{\{\log k-K-\Lambda-\log W\geq0\}}{\mathrm e}
^{-({1}/{\lambda)}(\log k-K-\Lambda-\log W)} \bigr)
\\
& \leq &\E \biggl( \mathbb{P}\biggl(E<\frac{{\mathrm e}^{K}}{k}W \Big\vert W\biggr) \biggr)+\E
\bigl( {\mathrm e}^{-({1}/{\lambda})(\log k-K-\Lambda-\log W)} \bigr)
\nonumber
\\
& \leq& O(1/k)\E(W)+O\bigl(k^{-1/\lambda}\bigr)\E \bigl( E^{-1/\lambda}W^{1/\lambda}
\bigr) =O\bigl(k^{-1/\lambda}\bigr)\nonumber,
\end{eqnarray}
and it follows that $k-\widehat{D}_k=O_{\prob}(k^{1-1/\lambda})$ as
in the
previous case.

To show the corresponding lower bound, let $\varepsilon>0$ be given,
and choose $K<\infty$ large enough that $\prob(M_k>\log
k+K)<\varepsilon
$. Similar to \eqref{eqDstarBinomialLower},
%
\begin{equation}
\label{eqDstarBinomialUpper}\qquad k-\widehat{D}_k \geq-1+\sum
_{i=1}^k \mathbh{1}_{\{\Lambda_i+\log
W_i+\lambda Y_i \geq\log k+K\}}\qquad \mbox{on }\{M
\leq\log k+K\}.
\end{equation}
We estimate
%
\begin{eqnarray}
&&\prob(\Lambda+\log W+\lambda Y \geq\log k+K)
\nonumber
\\[-8pt]
\\[-8pt]
\nonumber
&&\qquad \geq \prob( \Lambda\geq0)\prob(W\geq\delta)\prob \bigl( \lambda Y\geq \log k+K+
\log(1/\delta) \bigr) \geq c k^{-1/\lambda}
\end{eqnarray}
provided $\delta>0$ is small enough that $\prob(W\geq\delta)>0$.
Therefore, apart from an event of small probability, $k-\widehat
{D}_k+1$ is stochastically larger than a $\operatorname{Binomial}(k,ck^{-1/\lambda})$
random variable, and such a variable is itself $\Theta_\prob
(k^{1-1/\lambda})$.

The proof of part (d) is similar. For the upper
bound, it suffices to show that $1-p_k=O(k^{-1}\log k)$. Recall that
$\Lambda=-\log E$, and write the standard exponential variable $Y$ as
$Y=-\log U$, where $U$ is $\operatorname{Uniform}[0,1]$. Then
%
\begin{equation}\qquad
1-p_k = \prob(-\log E+\log W-\log U \geq\log k - K) = \prob\bigl(EU
\leq W{\mathrm e}^K/k\bigr).
\end{equation}
Splitting according to the value of $U$, we can then estimate
%
\begin{eqnarray}
\prob(EU\leq z)&\leq& z+\prob(U\geq z, E\leq z/U)
\nonumber
\\[-8pt]
\\[-8pt]
\nonumber
&\leq& z+\int_z^1
(z/u)\,\dd u = z\bigl(1+\log(1/z)\bigr),
\end{eqnarray}
so that $1-p_k \leq\E ( (W{\mathrm e}^K/k)(1+\log k-\log W -K)
 )$.
Note that the term $-W\log W$ is bounded above, so we conclude that
$1-p_k \leq O(k^{-1}\log k)$, as required. Similarly, for the lower
bound, we use $\P(E\leq y)\geq cy$ for $y\leq1$ to estimate $\prob
(EU\leq z)\geq\int_z^1 (cz/u)\,\dd u=cz\log(1/z)$ for any $z<1$, and we
conclude that
%
\begin{eqnarray}
\prob(\Lambda+\log W+Y\geq\log k+K) &\geq& \prob(W\geq\delta)\prob\bigl(EU\leq
\delta{\mathrm e}^K/k\bigr)
\nonumber
\\[-8pt]
\\[-8pt]
\nonumber
&\geq &c k^{-1} \log k
\end{eqnarray}
provided that $\prob(W\geq\delta)>0$ and that $k$ is large enough.
\end{pf*}

\subsection{Degree asymptotics: CM with infinite-variance degrees}

Now we prove that if the degrees in the configuration model have
infinite variance, then the shortest path tree reveals an asymptotic
proportion $p$ of the original degree.
The proof of Theorem~\ref{tCMdegreekinf} is similar to the proof of
Theorem~\ref{tCMdegreekfin}, except that here the asymptotic
proportion of revealed edges is $p<1$ and we need both upper and lower bounds.

\begin{pf*}{Proof of Theorem~\ref{tCMdegreekinf}}
Recall the notation $\xi_k=\min_{{i\in[k]}} (V_i +E_i)$. The
hypotheses on $V$ and $E$ imply that $\xi_k\stackrel{{\mathbb
P}}{\longrightarrow}0$ as $k\to\infty$.
Let $\varepsilon>0$ be given. Since $V$ and $E$ have continuous
distributions, we may choose $x>0$ such that $p-\varepsilon\leq\prob
(V-E>x)\leq\prob(V-E>0)=p$. Then
%
\begin{equation}
\label{eqDstarBinomialInfVar} \sum_{i=1}^k
\mathbh{1}_{\{V_i-E_i>x\}} \leq\widehat{D}_k \leq 1+\sum
_{i=1}^k \mathbh{1}_{\{V_i-E_i>0\}} \qquad\mbox{on }\{
\xi_k < x\},
\end{equation}
and each sum on the left-hand side of \eqref{eqDstarBinomialInfVar} is
$\operatorname{Binomial}(k,q)$ for some parameter $q\in[p-\varepsilon,p]$. Since
$\prob
(\xi_k < x)\to1$, it follows from the concentration of the binomial
distribution that $\prob(k(p-2\varepsilon)\leq\widehat{D}_k \leq
k(p+\varepsilon))\to1$, and since $\varepsilon>0$ was arbitrary, this
shows that $\widehat{D}_k=p \cdot k \cdot(1+o_{\prob}(1))$.
\end{pf*}

\subsection{Degree asymptotics: The complete graph}

In this section we prove Theorem~\ref{tKnDegreeTails}. This theorem
shows that the degree distribution on the shortest path tree $\TT_n$
behaves very differently for the complete graph $K_n$ compared to the
configuration model ${\mathrm{CM}}_n(\mathbf{d})$.

\begin{figure}

\includegraphics{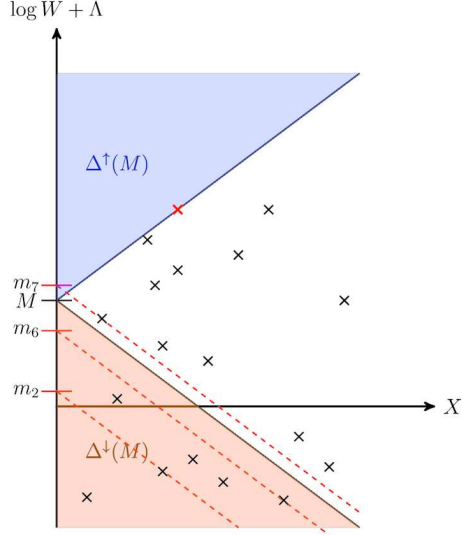}

\caption{The Poisson point process $\cP$. Crosses denote the points
$(X_i, \log W_i + \Lambda_i)$, and the coloured areas indicate the
upward- and downward-facing infinite triangles $\Delta^{\uparrow}(M)$
and $\Delta^{\downarrow}(M)$. The maximum $M$ of $\Lambda+ \log W +
\lambda_s X$ is taken at the thick red cross. By \protect\eqref{eqhatDandP},
the degree in this configuration is $1+\cP(\Delta^
\downarrow(M))=1+6=7$. The dashed lines indicate the values
$m_2,m_6,m_7$ introduced in the proof of
Theorem~\protect\ref{tKnDegreeTails}\textup{(b)} and \textup{(c)}.}\label{figmeasures}
\end{figure}

We use the representation of the limiting degree distribution from
Theorem~\ref{tKnDegrees}(a).
Recall that the points $ (X_i )_{i=1}^{k}$ form a Poisson point
process (PPP) with intensity measure $\mu_s(\dd x)=\frac
{1}{s}x^{1/s-1}\,\dd x$ on $\R^+$. Since $\Lambda_i+\log W_i$, $i\in\N$,
are i.i.d. random variables, the points $(X_i, \Lambda_i+\log W_i)$
form a PPP $\cP$ on $\R^+\times\R$ (see, e.g., \cite{Resnick86},
Proposition~2.2) with the product intensity measure
$\widetilde{\mu}_s$ {given by}
%
\begin{equation}
\label{eqintensitymeasure} \widetilde\mu_s(\dd x \,\dd y)=\mu_s(\dd
x)\cdot\Pv (\Lambda +\log W\in\dd y ).
\end{equation}
Let $\cP(S)$ stand for the number of points $(X_i, \Lambda_i+\log W_i)$
in this Poisson point process for any measurable set $S\subset\R
^+\times\R$.
We introduce infinite upward- and downward-facing triangles (see Figure~\ref{figmeasures}) with $y$-intercept $m$,
%
\begin{eqnarray}
\label{deftriangles} %
\Delta^{\uparrow}(m)& =& \bigl
\{(x,y) \in\R^+ \times\R \dvtx y\ge m+\lambda_s x \bigr\},
\nonumber
\\[-8pt]
\\[-8pt]
\nonumber
\Delta^{\downarrow}(m)& = &\bigl\{(x,y) \in\R^+ \times\R\dvtx y \le m-
\lambda_s x \bigr\}.
\end{eqnarray}
With this notation in mind, we can rewrite $M=\max_{i\in\N}(\La_i +
\log W_i - \lambda_s X_i)$
from~\eqref{eqKndegreecharacterisation} as $M=\sup\{m\in\R\dvtx \cP(\Delta^\uparrow(m))\geq1\}=\inf\{m\in\R\dvtx \cP(\Delta
^\uparrow(m))=0\}$, and
%
\begin{equation}
\label{eqMidentity} \Pv ( M \ge m ) = \Pv\bigl(\cP\bigl( \Delta^{\uparrow}(m)
\bigr)\ge1\bigr) =1- \exp\bigl\{ -\widetilde\mu_s\bigl(
\Delta^{\uparrow}(m)\bigr)\bigr\}.
\end{equation}

Thus, \eqref{eqKndegreecharacterisation} implies that
%
\begin{equation}
\label{eqhatDandP} \widehat D-1 = \sum_{i \in\N}
\mathbh{1}_{\{\lambda_s X_i + \log
W_i + \Lambda_i<M\}} = \cP\bigl(\Delta^\downarrow(M)\bigr).
\end{equation}
Moreover, by the Poisson property, conditional on $M$, the number $\cP
(\Delta^\downarrow(M))$ is Poisson with parameter $\widetilde\mu_s(
\Delta^{\downarrow}(M) )$ [since $\{M\geq m\}$ is measurable with
respect to the restriction of $\cP$ to $\Delta^\uparrow(m)$, whereas
$\Delta^\uparrow(m)
\cap
\Delta^\downarrow(m)=\{(0,m)\}$ has measure zero]. Hence, by the law
of total
probability,
%
\begin{equation}
\label{eqLTP} \Pv ( \widehat D -1 =k ) = \int_{-\infty}^{\infty}
\Pv \bigl( \Poi\bigl(\widetilde\mu_s\bigl( \Delta^{\downarrow}(m)
\bigr) \bigr) = k \bigr) \Pv (M\in\dd m ).
\end{equation}
Thus, in order to understand $\widehat D-1$, we need to investigate the
behavior of $\widetilde\mu_s( \Delta^{\downarrow}(m))$ and
$\widetilde
\mu_s( \Delta^{\uparrow}(m))$ as functions of $m$. We start with $s=1$,
in which case~\eqref{eqLTP} leads to analytically tractable integrals.

\begin{pf*}{Proof of Theorem~\ref{tKnDegreeTails}\normalfont{(a)}}
In this case, the weights are exponential, and the evolution of the
shortest path tree is the same as that of the \emph{Yule process}, and
\mbox{$W\stackrel{d}{=}E$}. Thus $-\log W -\Lambda\stackrel{d}{=}\Lambda
'-\Lambda$, with
$\Lambda', \Lambda$ i.i.d. Gumbel random variables.
The distribution of $\Lambda'-\Lambda$ is called the logistic
distribution and is clearly symmetric about $0$.
We compute
%
\begin{eqnarray}
\label{eqPDFofL} \Pv\bigl(\Lambda'-
\Lambda\ge x\bigr) &=& \Ev \bigl[\mathbb{P}\bigl(\Lambda\leq-x+
\Lambda' \vert\Lambda'\bigr) \bigr] = \Ev \bigl[ \exp
\bigl\{-{\mathrm e}^{x-\Lambda'}\bigr\} \bigr]
\nonumber
\\[-8pt]
\\[-8pt]
\nonumber
& =& \Ev \bigl[ \exp\bigl\{ -{\mathrm e}^{x} E \bigr\} \bigr]=
\frac{1}{1+{\mathrm e}^{x}}.
\end{eqnarray}
We have $\lambda_s=\lambda_1=1$ and $\mu_1(\mathrm {d}x)=\mathrm{ d} x$, so
%
\begin{eqnarray}
\label{eqmutri} \widetilde\mu_1 \bigl(\Delta ^{\uparrow}(m)
\bigr) &=& \int_{0}^{\infty} \Pv\bigl(
\Lambda'-\Lambda\ge x+m\bigr) \,\mathrm {d} x =\int_{0}^{\infty}
\frac{{\mathrm
e}^{-(m+x)}}{1+{\mathrm e}^{-(m+x)}} \,\mathrm {d} x
\nonumber
\\[-8pt]
\\[-8pt]
\nonumber
&=& \log \bigl(1+{\mathrm e}^{-m} \bigr).
\end{eqnarray}
Thus the distribution of $M$ is the same as that of $\La'-\La$,
%
\begin{equation}
\label{eqdistributionofmin} \Pv(M\ge m)= \Pv\bigl[ \cP\bigl( \Delta^{\uparrow}(m) \bigr)
\ge1\bigr] =1-{\mathrm e}^{- \widetilde\mu_1 ( \Delta^{\uparrow}(m)
)} =\frac
{1}{1+{\mathrm e}^m}.
\end{equation}
[In general, recall from Lemma~\ref{lemMLambdalogW} that $M\stackrel{d}{=}
\Lambda+\log W$; thus \eqref{eqdistributionofmin} is an expression of
the symmetry of $\Lambda'-\Lambda$ that is particular to the case
$s=1$.] Similarly,
%
\begin{eqnarray}
\label{eqmus1} \widetilde\mu_1 \bigl(\Delta^{\downarrow}(m) \bigr)
&= &\int_0^\infty\Pv\bigl(\La'-\La
\leq m-x\bigr)\, \mathrm {d}x
\nonumber
\\[-8pt]
\\[-8pt]
\nonumber
&= &\int_0^\infty
\frac{{\mathrm e}^{m-x}}{1+{\mathrm e}^{m-x}} \,\mathrm {d}x =
\log\bigl(1+{\mathrm e}^m\bigr).
\end{eqnarray}
Combining \eqref{eqLTP}, \eqref{eqdistributionofmin} and \eqref{eqmus1},
\begin{eqnarray*}
\Pv( \widehat{D}-1 = k) &=& \int_{-\infty}^\infty
\Pv \bigl(\Poi\bigl(\log \bigl(1+{\mathrm e}^m\bigr)\bigr)=k \bigr)
\,\mathrm {d }\Pv(M\le m)
\\
&=& \int_{-\infty}^\infty\frac{1}{1+{\mathrm e}^m}
\frac{(\log
(1+{\mathrm e}^m))^k}{k!} \frac{{\mathrm e}^m}{(1+{\mathrm e}^m)^2}
\,\mathrm{ d}m
\\
&= &\int_0^\infty\frac{t^k}{k!} {\mathrm
e}^{-2t} \,\mathrm{ d}t = \frac
{1}{2^{k+1}},
\end{eqnarray*}
where in the last line we used the change of variables $t = \log
(1+{\mathrm e}
^m)$. This completes the proof of Theorem~\ref{tKnDegreeTails}(a).
\end{pf*}

When $s\neq1$, we do not have a closed form for the distribution of
$W$, so we need to estimate the parameters of the Poisson variables in
\eqref{eqLTP}. The following lemma summarizes the asymptotic
properties of $\widetilde{\mu}_s( \Delta^{\downarrow}(m))$,
$\widetilde{\mu}_s( \Delta
^{\uparrow
}(m))$ and $M$ that we will need. To state it, we define $g$ to be the
inverse of the function $m\mapsto\widetilde{\mu}_s(\Delta
^\downarrow(m))$, and set $\delta
=\widetilde{\mu}_s
(\Delta^\downarrow(1))>0$.

%
\begin{Lemma}\label{lemtildemuBounds}
Fix $s>0$. Then:
\begin{longlist}[(a)]
\item[(a)]
uniformly over $m\geq1$ and $u\geq\delta$,
%
\begin{eqnarray}
\widetilde{\mu}_s\bigl(\Delta^\downarrow(m)\bigr) &=& (m/
\lambda_s)^{1/s}\bigl(1+O(1/m)\bigr), \label{eqtildemuDownAsymp}
\\
g(u) &=& \lambda_s u^s + O(1); \label{eqgAsymp}
\end{eqnarray}

\item[(b)]
there is a constant $c$ (depending on $s$) such that, for any $m\geq0$,
%
\begin{equation}
\label{eqtildemuUpBounds} c {\mathrm e}^{-m}\leq\widetilde{\mu}_s
\bigl(\Delta^\uparrow(m)\bigr) \leq {\mathrm e}^{-m}.
\end{equation}
Furthermore the random variable $M$ has a density $\frac{\prob(M\in
\dd
m)}{\dd m}$ with respect to Lebesgue measure, and
%
\begin{equation}
\label{eqMdensity} {\mathrm e}^{-m}\E\bigl(W {\mathrm e}^{-W}
\bigr) \leq\frac{\prob(M\in\dd
m)}{\dd m} \leq{\mathrm e}^{-m}, \qquad m\geq0;
\end{equation}

\item[(c)]
there is a constant $C$ (depending on $s$) such that, for any $m\geq1$,
%
\begin{equation}
\frac{\dd}{\dd m}\widetilde{\mu}_s\bigl(
\Delta^\downarrow(m)\bigr) \leq C m^{1/s}.
\end{equation}
\end{longlist}
\end{Lemma}

\begin{pf}
By the definition of $\widetilde{\mu}_s$,
%
\begin{eqnarray}
\label{eqtridownequality} \widetilde{\mu}_s\bigl(\Delta ^{\downarrow}(m)
\bigr) &=& \int_0^\infty\Pv( \Lambda+ \log W < m-
\lambda _s x) \,\dd\mu_s(x)
\nonumber
\\
&=&\int_0^\infty\E \bigl( \exp \bigl( -{\mathrm
e}^{-m+\lambda_s x}W \bigr) \bigr) \,\dd\mu_s(x)\\
& =&\int
_0^\infty\phi_{W}\bigl({\mathrm
e}^{-m+\lambda_s x}\bigr) \,\dd\mu_s(x),\nonumber
\end{eqnarray}
where $\phi_{W}(u)=\E({\mathrm e}^{-uW})$. We split the integral
into two
terms, and use the trivial bound $\phi_{W}(u)\leq1$ in the first
term to get
%
\begin{equation}
\label{eqtridownbound1} \widetilde\mu_s\bigl(\Delta^{\downarrow}(m)\bigr)
\le\int_0^{m/\lambda_s} 1 \cdot\dd\mu_s(x) +
\int_{m/\lambda_s}^\infty\phi_{
W}\bigl({\mathrm
e}^{\lambda
_s x-m}\bigr) \,\dd\mu_s(x).
\end{equation}
The first term equals $(m/\lambda_s)^{1/s}$, so we continue by showing
that the second term in \eqref{eqtridownbound1} is of smaller order.
Recall that $\phi_{W}$ satisfies the recursive relation \eqref
{eqphiWrecursion}.
By the monotonicity property of $\phi_{W}$, we have $\phi_{
W}(u{\mathrm e}^{-\lambda_s x})\leq\phi_{W}(1)$ as long as
$x\leq(\log
u)/\lambda_s$. Hence, for $u\geq1$,
%
\begin{eqnarray}
\label{eqmomgenbound} \phi_{W}(u) &\le&\exp \biggl\{ - \int
_{0}^{(\log u)/\lambda_s} \bigl(1-\phi_{W}(1)\bigr)
\,\dd\mu_s(x) \biggr\}
\nonumber
\\[-8pt]
\\[-8pt]
\nonumber
 &= &\exp \biggl\{ - \frac{(1-\phi_{W}(1)) (\log u)^{1/s}}{ \lambda
_s^{1/s}} \biggr\}.
\end{eqnarray}
Recalling the definition \eqref{eqmusDefinition} of $\mu_s$ and making
the subsitution $t=x-m/\lambda_s$, we conclude that the second term of
\eqref{eqtridownbound1} is at most
%
\begin{equation}
\int_0^\infty\exp \bigl( -\bigl(1-
\phi_{W}(1)\bigr)t^{1/s} \bigr) \frac
{1}{s} \biggl( t+
\frac{m}{\lambda_s} \biggr)^{1/s -1}\,\dd t.
\end{equation}
For $s>1$, the estimate $(t+m/\lambda_s)^{1/s-1}\leq(m/\lambda
_s)^{1/s-1}$ shows that the second term of \eqref{eqtridownbound1} is
$O(m^{1/s-1})$. For $s<1$, the bound $(t+m/\lambda_s)^{1/s-1}\leq
(2t)^{1/s-1} + (2m/\lambda_s)^{1/s-1}$ shows that the second term of
\eqref{eqtridownbound1} is $O(1)+O(m^{1/s-1})$, which is
$O(m^{1/s-1})$ uniformly over $m\geq1$. In either case we have
verified \eqref{eqtildemuDownAsymp}. By the definition of $\delta$,
\eqref{eqgAsymp} follows from \eqref{eqtildemuDownAsymp}, and this
proves part (a).

For part (b), the upper bound in \eqref
{eqtildemuUpBounds} follows from the bounds $\prob(\Lambda\geq
x)=1-{\mathrm e}
^{-{\mathrm e}^{-x}}\leq{\mathrm e}^{-x}$ {as follows:}
%
\begin{eqnarray}\label{eqtildemuUpFormula}
\widetilde{\mu}_s\bigl(\Delta^\uparrow(m)\bigr) &=& \int
_0^\infty\prob(\Lambda+\log W\geq m+
\lambda_s x) \,\dd\mu_s(x)
\nonumber
\\[-8pt]
\\[-8pt]
\nonumber
& \leq& \int_0^\infty\E\bigl({\mathrm
e}^{-m-\lambda_s x+\log W}\bigr) \,\dd\mu_s(x) = {\mathrm e}^{-m},
\nonumber
\end{eqnarray}
since $\E(W)=1=\int_0^\infty{\mathrm e}^{-\lambda_s x}\,\dd\mu
_s(x)$. For the
lower bound, note that $\prob(W\geq1)>0$ [since $\E(W)=1$], so the
bound $\prob(\Lambda\geq x)\geq c'{\mathrm e}^{-x}$ gives
%
\begin{equation}
\widetilde{\mu}_s\bigl(\Delta^\uparrow(m)\bigr)\geq
\int_0^1 \prob(W\geq 1)\prob(\Lambda\geq m+
\lambda_s)\,\dd \mu_s(x)\geq c {\mathrm e}^{-m}.
\end{equation}
For \eqref{eqMdensity}, use Lemma~\ref{lemMLambdalogW} to express the
density of $M$ in terms of the density ${\mathrm e}^{-{\mathrm e}^{-x}}
{\mathrm e}^{-x}\,\dd x$ of
a Gumbel random variable,
%
\begin{eqnarray}
\prob(M\in\dd m) &=& \E \bigl( \mathbb{P}(\Lambda+\log W\in\dd m \vert W) \bigr)
\nonumber
\\[-8pt]
\\[-8pt]
\nonumber
& =& \E \bigl( {\mathrm e}^{-{\mathrm e}^{-m+\log W}} {\mathrm e}^{-m+\log
W} \bigr) \,\dd m
= \E \bigl( W{\mathrm e}^{-W{\mathrm e}^{-m}} \bigr) {\mathrm e}^{-m} \,\dd m.\hspace*{-20pt}
\end{eqnarray}
We may then bound $\E(W{\mathrm e}^{-W{\mathrm e}^{-m}})$ above and
below by $\E(W)=1$
and $\E(W{\mathrm e}^{-W})$, respectively, completing the proof of
\eqref
{eqMdensity} and part (b).

Finally, for part (c), note from \eqref
{eqtridownequality} that
%
\begin{equation}
\frac{\mathrm d}{\mathrm {d}m} \widetilde{\mu}_s\bigl(\Delta^\downarrow (m)
\bigr) = \int_0^\infty {\mathrm e}^{\lambda
_s x-m}
\bigl(-\phi_{W}'\bigl({\mathrm e}^{\lambda_s x-m}\bigr)
\bigr) \,\dd\mu_s(x).
\end{equation}
Recalling \eqref{eqphiWrecursion} and using the trivial bound $-\phi
_{W}'(u)\leq-\phi_{W}'(0)=\E(W)=1$,
%
\begin{equation}\qquad
\label{eqphiphi} \frac{\phi_{W}'(u)}{\phi_{W}(u)} = \int_0^\infty{
\mathrm e}^{-\lambda
_s x}\bigl(-\phi_{W}'\bigl(u{\mathrm
e}^{-\lambda_s x}\bigr)\bigr)\,\dd\mu_s(x) \leq \int
_0^\infty{\mathrm e}^{-\lambda_s x}\,\dd
\mu_s(x)=1,
\end{equation}
and using \eqref{eqmomgenbound} we conclude that
%
\begin{eqnarray}
\label{eqmuDerivativeBound} &&\frac{\dd}{\dd m} \widetilde{\mu}_s\bigl(
\Delta^\downarrow(m)\bigr)\nonumber\\
 &&\qquad\leq \int_0^\infty{
\mathrm e}^{\lambda_s x-m}\phi_{W}\bigl({\mathrm e}^{\lambda_s x-m}
\bigr) \,\dd\mu_s(x)
\nonumber
\\[-8pt]
\\[-8pt]
\nonumber
&&\qquad \leq \mu_s[0,m/\lambda_s]+\int_{m/\lambda_s}^\infty{
\mathrm e}^{\lambda
_s x-m}\exp \bigl(-c(\lambda_s x-m)^{1/s}
\bigr) \,\dd\mu_s(x)
\\
&&\qquad = (m/\lambda_s)^{1/s} + \int_0^\infty{
\mathrm e}^{z-cz^{1/s}}\frac
{(z+m)^{1/s-1}}{s\lambda_s^{1/s-1}}\,\dd z,\nonumber
\end{eqnarray}
where $z=\lambda_s x-m$. As before, we either bound $(z+m)^{1/s-1}\leq
m^{1/s-1}$ (if $s>1$) or $(z+m)^{1/s-1}\leq(2z)^{1/s-1}+(2m)^{1/s-1}$
(if $s<1$) to conclude that the last term in \eqref
{eqmuDerivativeBound} is $O(m^{1/s-1})+O(1)$. Hence the upper bound in
\eqref{eqmuDerivativeBound} is $O(m^{1/s})$ uniformly over $m\geq1$,
which completes the proof.
\end{pf}

With Lemma~\ref{lemtildemuBounds} in hand, we can now prove Theorem~\ref{tKnDegreeTails}(b) and (c).

\begin{pf*}{Proof of Theorem~\ref{tKnDegreeTails}\normalfont{(b)} and
\normalfont{(c)}}
From \eqref{eqLTP}, we see that the unlikely event $\{ \smash
{\widehat{D}-1=k}\}$ is achieved when the variables $M$ or $\Poi
(\widetilde{\mu}_s
(\Delta^\downarrow(m)))$, or both, are unusually large. As a
heuristic to evaluate
the costs of these alternatives, we can use Lemma~\ref
{lemtildemuBounds}(a) and (b) to approximate $\widetilde{\mu}_s(\Delta
^\uparrow(m))\approx{\mathrm e}^{-m}$,
$\widetilde{\mu}_s(\Delta^\downarrow(m))\approx(m/\lambda
_s)^{1/s}\mathbh{1}_{\{m\geq0\}}$,
leading to
%
\begin{eqnarray}
\label{DhatHeuristic} \prob(\widehat{D}-1=k) &\approx &\int_0^\infty
\frac{{\mathrm e}^{-(m/\lambda_s)^{1/s}}(m/\lambda
_s)^{k/s}}{k!}\bigl({\mathrm e} ^{-m} \,\dd m\bigr)
\nonumber
\\[-8pt]
\\[-8pt]
\nonumber
& =& \int_0^\infty\frac{s\lambda_s u^{s-1}}{k!} \exp \bigl(
-u-\lambda _s u^s +k\log u \bigr) \,\dd u
\end{eqnarray}
after the substitution $u=(m/\lambda_s)^{1/s}$. The exponential in
\eqref{DhatHeuristic} is maximized when $u=u_*$, where $u_*$ is the
unique solution of
%
\begin{equation}
u_*+s\lambda_s u_*^s = k.
\end{equation}
For $s<1$, we have $u_*\approx k$, corresponding to $m_*\approx\lambda
_s k^s$, whereas for $s>1$ we have $u_*\approx(k/s\lambda_s)^{1/s}$,
corresponding to $m_*\approx k/s$.

We now formalize this heuristic argument. For $k\in\N$, define the
random variables
%
\begin{equation}
m_k = \inf\bigl\{m\in\R\dvtx \cP\bigl(
\Delta^\downarrow(m)\bigr)\geq k\bigr\}.
\end{equation}
(See Figure~\ref{figmeasures}: $m_k$ is the value on the vertical axes
where the $k$th point enters the downward-facing triangle). Note that
each $m_k$ is a stopping time with respect to the filtration $(\sigma
(\cP\vert_{\Delta^\downarrow(m)}))_{m\in\R}$ generated by the
restrictions of
$\cP$
to $\Delta^\downarrow(m)$, $m\in\R$. In terms of $m_k$, we have
%
\begin{equation}
\label{eqDboundIffMbound} \{\widehat{D}-1\geq k\} = \{M\geq m_k\}.
\end{equation}
Since $\Delta^\uparrow(m)$ is disjoint from $\Delta^\downarrow(m)$,
it follows that
%
\begin{eqnarray}\label{eqMgeqmk}
 \mathbb{P}(M\geq m_k \vert
m_k=m)&=&\prob(M\geq m)=\prob\bigl(\cP\bigl(\Delta ^\uparrow(m)
\bigr)>0\bigr)
\nonumber
\\[-8pt]
\\[-8pt]
\nonumber
&=&1-{\mathrm e} ^{-\widetilde{\mu}_s(\Delta^\uparrow(m))}\leq\widetilde{\mu }_s\bigl(
\Delta^\uparrow(m)\bigr).
\end{eqnarray}

Since the function $m\mapsto\widetilde{\mu}_s(\Delta^\downarrow
(m))$ is continuous, the
sequence\break  $(\widetilde{\mu}_s(\Delta^\downarrow(m_k)))_{k=1}^\infty
$ forms a Poisson point
process on $(0,\infty)$ of intensity $1$. [This fact, which is
elementary to verify, is the analogue of the statement that applying a
continuous distribution function to a variable having that distribution
gives a $\operatorname{Uniform}(0,1)$ random variable.] In particular, $\widetilde
{\mu}_s(\Delta^\downarrow
(m_k))$ has the Gamma$(k,1)$ distribution with density $\Gamma(k)^{-1}
u^{k-1} {\mathrm e}^{-u}\,\dd u$.

For the upper bound, it suffices to estimate $\prob(\widehat{D}-1\geq
k)$. By \eqref{eqDboundIffMbound}, this amounts to bounding $\prob
(M\geq m_k)$. We begin with $s<1$, in which case the above heuristics
suggest that the dominant contribution to $\prob(\widehat{D}-1\geq k)$
comes when $\widetilde{\mu}_s(\Delta^\downarrow(m_k))\approx k$.
Partitioning according to the
value $u=\widetilde{\mu}_s(\Delta^\downarrow(m_k))$, and combining
with the fact that $\widetilde{\mu}_s
(\Delta^\downarrow(m_k))$ has the Gamma distribution, we obtain
%
\begin{eqnarray}
\prob(\widehat{D}-1\geq k) &\leq &\prob\biggl(\widetilde{\mu}_s\bigl(
\Delta^\downarrow(m_k)\bigr)\notin\biggl[\frac
{1}{2}k,
\frac{3}{2}k\biggr]\biggr)
\nonumber
\\
&&{} + \prob\biggl(\widetilde{\mu}_s\bigl(\Delta^\downarrow(m_k)
\bigr)\in\biggl[\frac
{1}{2}k,\frac{3}{2}k\biggr], M\geq g\bigl(
\widetilde{\mu}_s\bigl(\Delta^\downarrow(m_k)\bigr)
\bigr)\biggr)
\nonumber
\\[-8pt]
\\[-8pt]
\nonumber
& =& \prob\biggl(\mbox{Gamma}(k,1)\notin\biggl[\frac{1}{2}k,
\frac{3}{2}k\biggr]\biggr) \\
&&{}+ \int_{k/2}^{3k/2}
\frac{u^{k-1} {\mathrm e}^{-u}}{(k-1)!} \prob\bigl(M\geq g(u)\bigr)
\,\dd u,\nonumber
\end{eqnarray}
where we used that $g$ is the inverse function of $m\mapsto\widetilde
\mu_s(\Delta^\downarrow(m))$. We can continue estimating the
right-hand side as
%
\begin{eqnarray}\qquad
\label{eqDhatUpperBoundInitial} \prob(\widehat{D}-1\geq k)&\leq &{\mathrm e}^{-ck} + \int
_{k/2}^{3k/2} \frac{u^{k-1}{\mathrm
e}^{-u}}{(k-1)!} \widetilde{
\mu}_s\bigl(\Delta ^\uparrow\bigl(g(u)\bigr)\bigr) \,\dd u
\nonumber
\\[-8pt]
\\[-8pt]
\nonumber
& \leq& {\mathrm e}^{-ck}+\int_{k/2}^{3k/2}
\frac{\exp ( (k-1)\log u
-u-\lambda_s
u^s +O(1)  )}{(k-1)!} \,\dd u,
\end{eqnarray}
where we used that $\widetilde{\mu}_s(\Delta^\uparrow(g(u)))\le
\mathrm e^{-g(u)}$
by \eqref{eqtildemuUpBounds} and then the bound on $g(u)$ in \eqref
{eqgAsymp}.

Uniformly over the range of integration, Stirling's approximation and a
Taylor expansion give
\[
(k-1)\log u-u-\log\bigl((k-1)!\bigr)\leq-\frac{1}{8k}(k-1-u)^2+O(
\log k),
\]
whereas $\lambda_s u^s=\lambda_s (k-1)^s+O( (k^{s-1})(k-1-u))$. Hence
%
\begin{eqnarray}
\label{eqDhatUpperBound}&& \prob(\widehat{D}-1\geq k)
\nonumber
\\
&&\qquad
\leq {\mathrm e}^{-ck}
+{\mathrm e}^{-\lambda_s (k-1)^s + O(\log k)} \int_{k/2}^{3k/2} \exp
\biggl(-\frac
{(k-1-u)^2}{8k}\\
&&\hspace*{198pt}{}+O \bigl(k^{s-1}\bigr)\vert k-1-u\vert \biggr)
\,\dd u.\nonumber
\end{eqnarray}
The integral in \eqref{eqDhatUpperBound} is $\exp(O(k^{2s-1}))$ (this
can be seen by maximising the integrand), which is negligible compared
to $\exp(-\lambda_s k^s)$ since $s<1$, and this proves the upper bound.

For $s>1$, the dominant contribution to $\prob(\widehat{D}-1=k)$ is
expected to come when $u=\widetilde{\mu}_s(\Delta^\downarrow(m_k))$
satisfies $u\approx
(k/s\lambda_s)^{1/s}\ll k$. We partition into the events $\{u\geq k\}
$ (in which case we must have $M\geq m_k=g(u)\geq g(k)$), $\{u\leq
\delta=\widetilde{\mu}_s(\Delta^\downarrow(1))\}$ (in which case
we must have $m_k\leq1$ and
$\cP(\Delta^\downarrow(1))\geq k$), and $\{\delta\leq u\leq k\}$.
As in~\eqref
{eqDhatUpperBoundInitial}--\eqref{eqDhatUpperBound},
%
\begin{eqnarray}
 \prob(\widehat{D}-1\geq k) &\leq& \prob\bigl(M\geq
g(k)\bigr) + \prob\bigl(\cP\bigl(\Delta^\downarrow(1)\bigr)\geq k\bigr)
\nonumber
\\
& &{}+ \int_\delta^k \frac{\exp \{ (k-1)\log u -u-\lambda_s u^s +O(1)
 \}}{(k-1)!} \,\dd u
\nonumber
\\
& \leq &\widetilde{\mu}_s\bigl(\Delta^\uparrow\bigl(g(k)\bigr)
\bigr)+\prob\bigl(\Poi(\delta)\geq k\bigr)
\nonumber
\\[-8pt]
\\[-8pt]
\nonumber
&&{}+O(k)\frac
{\exp
 \{ \max_{\delta\leq u\leq k} (k\log u-\lambda_s u^s)  \}}{(k-1)!}
\\
& \leq& {\mathrm e}^{-\lambda_s k^s+O(1)}+{\mathrm e}^{-k\log
k+O(k)}\nonumber\\
&&{}+O
\bigl(k^2\bigr)\frac{\exp\{
({k}/{s})\log({k}/{(s\lambda_s)})-{k}/{s}\}}{k!},\nonumber
\end{eqnarray}
where we used \eqref{eqMgeqmk} first and then \eqref{eqgAsymp} to
bound $g(k)$.
The desired bound follows by Stirling's approximation.

For the lower bound, let $\varepsilon>0$ be given. We begin with $s>1$.
By Lem\-ma~\ref{lemtildemuBounds}(a), uniformly
over $m\in[k,k^{1+\varepsilon}]$, we have $\widetilde{\mu}_s(\Delta
^\downarrow(m)
)=k^{1/s+O(\varepsilon)}$. Therefore, using~\eqref{eqLTP} and
Stirling's approximation,
%
\begin{eqnarray}
 \mathbb{P}(\widehat{D}-1=k \vert M=m) &=& \exp \bigl
\{k\log\widetilde {\mu}_s\bigl(\Delta^\downarrow(m) \bigr)-
\widetilde{\mu}_s \bigl(\Delta^\downarrow(m)\bigr) \bigr\}/k!
\nonumber
\\[-8pt]
\\[-8pt]
\nonumber
&=&\exp\bigl\{\bigl(1/s-1+O(\varepsilon)\bigr)k\log k\bigr\}.
\end{eqnarray}
On the other hand, to estimate $\Pv(M\in\mathrm[k, k^{1+\ve}])$ write
%
\begin{eqnarray}
&&\bigl\{k\leq M\leq k^{1+\varepsilon}\bigr\}
\nonumber
\\[-8pt]
\\[-8pt]
\nonumber
&&\qquad=\bigl\{\cP\bigl(\Delta^\uparrow
\bigl(k^{1+\varepsilon}\bigr)\bigr)=0\bigr\}\cap\bigl\{\cP\bigl(
\Delta^\uparrow(k)\setminus \Delta^\uparrow\bigl(k^{1+\varepsilon}
\bigr)\bigr)>0\bigr\}.
\end{eqnarray}
By Lemma~\ref{lemtildemuBounds}(b),
$\widetilde{\mu}_s
(\Delta
^\uparrow(k^{1+\varepsilon}))\leq{\mathrm e}^{-k^{1+\varepsilon
}}\to0$, so the
first event on the right-hand side occurs with high probability as
$k\to\infty$. Since in addition $\widetilde{\mu}_s(\Delta^\uparrow
(k))\geq c
{\mathrm e}
^{-k} \gg\widetilde{\mu}_s(\Delta^\uparrow(k^{1+\epsilon}))$, it
follows that the
second event occurs with probability at least $c' {\mathrm e}^{-k}$. Combining
all of these estimates gives the result.

Similarly, for $s<1$, let $m\in[g(k),g(k+1)]$, and set $u=\widetilde
{\mu}_s
(\Delta^\downarrow(m)
)$, so that $u\in[k,k+1]$. Uniformly over this range, we have $\log
u=\log k+o(k^s)$, and it follows using Stirling's approximation that
\[
\mathbb{P}(\widehat{D}-1=k \vert M=m)=\exp \{ -u+k\log u \}/k! = \exp \bigl\{ o
\bigl(k^s\bigr) \bigr\}.
\]
By Lemma~\ref{lemtildemuBounds}(b), we have
\[
\prob\bigl(g(k)\leq M\leq g(k+1)\bigr)\geq c {\mathrm e}^{-g(k+1)}
\bigl(g(k+1)-g(k)\bigr).
\]
We have $g(k+1)\sim g(k)\sim\lambda_s k^s$ by Lemma~\ref
{lemtildemuBounds}(a). To bound $g(k+1)-g(k)$,
note that the definition of $g$ implies
%
\begin{equation}
\label{eqgDifferenceBound} \bigl( g(k+1)-g(k) \bigr) \cdot\max_{g(k)\leq m\leq g(k+1)}
\frac
{\dd
}{\dd m}\widetilde{\mu}_s\bigl(\Delta^\downarrow(m)
\bigr) \geq1.
\end{equation}
We apply Lemma~\ref{lemtildemuBounds}(c)
with $m\sim\lambda_s k^s$, so that \eqref{eqgDifferenceBound} gives
$g(k+1)-g(k)\geq c/k$. Consequently $\prob(g(k)\leq M\leq g(k+1))\geq
{\mathrm e}
^{-\lambda_s k^s+o(k^s)}$, and this completes the proof.
\end{pf*}
%

\section{Deterministic edge weights}\label{sY=1Proof}

In this section we prove Theorem~\ref{tY=1}. The proof has some
similarity to the proofs in Section~\ref{spart-c-proofs}.

\begin{pf*}{Proof of Theorem~\ref{tY=1}}
Write $f(z)=\E(z^D)$ for the generating function of the degree
distribution $D$. It suffices to show that the generating function for
$\widehat{D}$ matches with the expression in \cite{clauset-newman08}, equation
(1),
%
\begin{equation}
\label{eqY=1ToProve} \E \bigl( z^{\widehat{D}} \bigr) = z\int_0^1
f' \biggl( t-(1-z)\frac
{f' ( {f'(t)}/{f'(1)}  )}{f'(1)} \biggr) \,\dd t.
\end{equation}
Since $Y=1$, we have ${\mathrm e}^{-\lambda Y}=1/\nu$, and recursive equation
\eqref{eqphiWrecursionCMFinVar} becomes
%
\begin{equation}
\label{eqphiWrecursionDeterministic} \phi_{W}(u)=\frac{f'(\phi_{W}(u/\nu))}{f'(1)}.
\end{equation}
Using symmetry, writing $\Lambda_1=-\log E_1$ and recalling that
$\prob
(\Lambda+\log W<x)=\phi_{W}({\mathrm e}^{-x})$,
%
\begin{eqnarray}
\E \bigl( z^{\widehat{D}} \bigr)
& = &\sum_{i=2}^\infty\prob(D=i)\sum
_{k=1}^i z^k i\pmatrix{i-1
\cr
k-1} \nonumber\\
&&\hspace*{68pt}{}\times\prob
\pmatrix{
M=\Lambda_1+\log W_1-\log\nu;
\vspace*{2pt}\cr
\Lambda_j+\log W_j+\log\nu<M\mbox{ for }j=2,\ldots,k;\mbox{ and}
\vspace*{2pt}\cr
\Lambda_j+\log W_j-\log\nu>M\mbox{ for }j=k+1,\ldots,i}\hspace*{-6pt}
%
\nonumber
\\
& =& \sum_{i=2}^\infty i\prob(D=i)\sum
_{k=1}^i z^k\pmatrix{i-1
\cr
k-1}\nonumber\\
&&\hspace*{58pt}{}\times \E \bigl( \mathbb{P}(\Lambda_j+\log W_j<
\Lambda_1+\log W_1-2\log\nu \vert\Lambda_1,W_1)^{k-1}
\nonumber\\
&&\hspace*{80pt}{}\times \mathbb{P}(\Lambda_1+\log W_1-2\log\nu<
\Lambda_j+\log W_j\\
&&\hspace*{191pt}{}<\Lambda _1+\log
W_1 \vert\Lambda_1,W_1)^{i-k}
\bigr)
\nonumber
\\
& =& z\E \Biggl( \sum_{i=2}^\infty i
\prob(D=i) \bigl( z\phi_{W}\bigl(\nu^2
E_1/W_1\bigr)\nonumber\\
&&\hspace*{83pt}{} + \bigl( \phi_{W}(E_1/W_1)
- \phi_{W}\bigl(\nu^2 E_1/W_1
\bigr) \bigr) \bigr)^{i-1} \Biggr)
\nonumber
\\
& =& z\E \bigl( f' \bigl( \phi_{W}(E_1/W_1)
- (1-z)\phi_{
W}\bigl(\nu ^2 E_1/W_1
\bigr) \bigr) \bigr).\nonumber
\end{eqnarray}
Applying \eqref{eqphiWrecursionDeterministic} twice, we obtain
%
\begin{equation}\qquad
\label{eqzhatDFormulaWithphiW} \E \bigl( z^{\widehat{D}} \bigr) = z\E \biggl( f'
\biggl( \phi_{W}(E_1/W_1) - (1-z)
\frac{f' (
{f'(\phi_{W}(E_1/W_1))}/{f'(1)}  )}{f'(1)} \biggr) \biggr).
\end{equation}
Finally, since $W$ is positive and finite-valued, $\phi_{
W}^{-1}(t)$ is defined for each $t\in(0,1)$, and we can compute
%
\begin{eqnarray}
\prob\bigl( \phi_{W}(E_1/W_1) < t
\bigr)& =& \prob\bigl( E_1 > W_1 \phi_{
W}^{-1}(t)
\bigr) = \E \bigl( {\mathrm e}^{-W_1 \phi_{W}^{-1}(t)} \bigr)
\nonumber
\\[-8pt]
\\[-8pt]
\nonumber
& =& \phi _{W}
\bigl(\phi_{W}^{-1}(t)\bigr)=t,
\end{eqnarray}
so that $\phi_{W}(E_1/W_1)$ has the $\operatorname{Uniform}(0,1)$ distribution.
Thus the expectation over the value of $\phi_W(E_1/W_1)$ in \eqref
{eqzhatDFormulaWithphiW} is equivalent to the integration in \eqref
{eqY=1ToProve}.
\end{pf*}

\section*{Acknowledgments}
S. Bhamidi would like to thank  Eurandom for their hospitality  where this work
commenced.\
We thank James Wilson and Frances Tong for their help in
creating the pictures in Figures~\ref{figsimulations} and \ref
{figsimul-random-reg}. We thank two anonymous referees whose
suggestions helped with the clarity of the paper.


%





\printaddresses
\end{document}